\documentstyle{amsppt}

\def\exam#1#2{\dimen1=3.25em
                      \noindent\hangindent=\dimen1\hangafter1
                      \hbox to\dimen1{#1\hfil~~}#2} 
\baselineskip=15pt plus 2pt
\magnification=1200
\def\cl #1{\cal{#1}}

\def\ov #1{\overline{#1}}

\def\bs{\backslash}
\def\n{\noindent}

\def\bx{$\square$}
\overfullrule=0pt

\def\cl #1{{\Cal{#1}}}
\def\card{\,\text{card}\,}
\def\diam{\,\text{diam}\,}

\def\conv{\text{conv}}
\def\lno{\leqalignno}
\def\l{\langle}
\def\r{\rangle}

\leftheadtext{Michel Talagrand}

\topmatter

\title Constructions of Majorizing Measures,\\
Bernoulli processes and cotype \endtitle

\rightheadtext {Majorizing Measures, Bernoulli processes and cotype}
\affil February 1994\endaffil

\author Michel Talagrand\endauthor

\affil University of Paris VI and The Ohio State University\endaffil
\address Equipe d Analyse-Tour 56, E.R.A. au C.N.R.S. no. 754, Universit\'e 
Paris VI, 4 Pl Jussieu, 75230  Paris Cedex 05, 
FRANCE  and  
Department of Mathematics, The Ohio State University, 231 W. 18th Ave., 
Columbus, OH  43210-1174 
\endaddress

\abstract We present three methods to construct majorizing measures in 
various settings.
 These methods are based on direct constructions of increasing sequences of
partitions through a simple exhaustion procedure rather than on the
construction of well separated ultrametric subspaces.  The first scheme of
construction provides a simple unified proof of the Majorizing 
Measure Theorem for Gaussian processes and of the following fact.  If $A,B$ 
are balanced  convex sets in a vector space, and if $A$ is sufficiently 
convex, a control of the covering numbers $N(A,\varepsilon B)$ for all 
$\varepsilon>0$ implies the (a priori stronger) existence of a majorizing 
measure on $A$ provided with the distance induced by $B$.  This establishes, 
apparently for the first time, a clear link between geometry and majorizing 
measures, and generalizes the  earlier results on majorizing measures on 
ellipsoids in Hilbert space, that were obtained by specific methods.  Much 
of the rest of the paper is concerned with the structure of bounded 
Bernoulli (=Radmacher) processes.  The main conjecture on their structure 
is reformulated in several ways, that are shown to be equivalent, and to be 
equivalent to the existence of certain majorizing measures.  Two schemes of 
construction of majorizing measures related to this problem are 
presented.  One allows to describe Bernoulli processes
when the index set, provided with the supremum norm, is sufficiently small. 
The other allows to prove a weak form of the main conjecture.  This result,
while not sufficient to characterize boundedness of Bernoulli processes,
 allows to prove the
remarkable fact that for any continuous operator $T$ from $C(K)$ to $E$, 
the Rademacher cotype-2 constant of $T$ is controlled by the maximum of the 
Gaussian cotype-2 constant of $T$ and of its $(2,1)$-summing norm.  It is 
also proved, as a consequence of one of the main inequalities on Bernoulli 
processes, that in a Banach space $E$ of dimension $n$, at most 
$n\log n \log\log n$ vectors suffices to compute the Rademacher cotype 
$2$ constant of $E$ within a universal constant. 
\endabstract 


\endtopmatter

\n{\bf 1 - Introduction 

The notion of majorizing measure has allowed considerable progress in the
study of certain stochastic processes, in particular Gaussian processes.

Given a metric space $(T,d)$, and a probability
measure $\mu$ on $T$, we set 
$$\lno{\gamma_{1/2}(T,d,\mu)=\sup\limits_{x\in T}
\int^{\infty}_0\sqrt {\log{1\over \mu
(B(x,\varepsilon))}}d\varepsilon
,&&(1.1)\cr}$$
where $B(x,\varepsilon)$ is the closed ball for $d$ centered at $x$ of 
radius $\varepsilon$.  One should observe that, since $\log 1=0$, the 
integrand is zero when $\varepsilon$ is larger than the diameter of $T$.  
We set 
$$\lno{\gamma_{1/2}(T,d) =\inf\limits_{\mu}\gamma_{1/2}
(T,d,\mu),&&(1.2)\cr}$$
where the infimum is taken over all possible choices of $\mu$.  It is
implicitly assumed in (1.1) that the ball $B(x,\varepsilon)$ is
$\mu$-measurable.  It must be pointed out that in (1.1) it is equivalent to
assume that $\mu$ is supported by a countable subset of $T$.  Actually, 
despite their name, majorizing measures have little connection with measure 
theory, and are actually a kind of weight system to measure the size of $T$.  
 The set $T$ will often be a subset of a Hilbert
space $H$, $d$ being the distance induced from $H$.  In that case we will 
write
$\gamma_{1/2}(T)$ rather than $\gamma_{1/2}(T,\Vert\cdot\Vert_2)$.

Consider an integer $M$, and standard independent Gaussian random variables
$(g_i)_{i\le M}$.  For a subset $T$ of $\Bbb R^M$, we set 
$$G(T)=\sup\{E\sup\limits_{t\in S}\sum_{i\le M}t_ig_i;\quad S\subset T;~~
S~~\text{finite}\}$$

\n where $t\in\Bbb R^M$ is written as $t=(t_i)_{i\le M}$.  It is proved in 
[T1] that for a certain constant $K$, we have 
$$\leqalignno{{1\over K}\gamma_{1/2}(T)\le G(T)\le K\gamma_{1/2}(T),
&&(1.3)\cr}$$

\n when $T$ is provided by the distance induced by $\ell^M_2$.  (The 
right-hand inequality is an earlier result of X. Fernique; the left hand 
side is known as the majorizing measure theorem).  We observe that (1.3) 
does not depend on $M$.  By approximation one could thus work in 
$\Bbb R^{\Bbb N}$.  But the present setting offers the advantage that one 
does not have to bother about infinite series.

A simpler proof of (1.3) is given in [T4].  An even simpler proof will be
given in Section 2 of the present paper.  While (1.3) is now fairly easy 
to prove, the construction of majorizing measures (i.e. of measures on $T$ 
that witness the left-hand side inequality of (1.3) in practical situations) 
is a difficult question.  One reason is that, while (1.3) is in principle 
(as explained in [T2]) a theorem about geometry of the 
Hilbert space, this geometric aspect is not understood.  There are 
situations (in particular some concrete classes of functions on $[0,1]^2$ 
that are studied in [T7]) where one has a rather geometrical knowledge of 
$T$ and where the precise computation of $\gamma_{1/2}(T)$ is currently 
intractable.  The first main contribution of the present paper will be the 
description in a simple but important situation of a precise link between 
the geometry of $T$ and the value of certain functionals 
$\gamma_{\alpha, \beta} (T)$ that generalize $\gamma_{1/2}(T)$ and that 
have been introduced in [T7].  Given a metric space $(T,d)$, numbers 
$\alpha, \beta >0$, and a probability measure $\mu$ on $T$, we set
$$\lno {\gamma_{\alpha, \beta}(T,d, \mu) = \sup\limits_{x \in T} 
(\int^\infty_0 \varepsilon^\beta \left( \log {1 \over \mu (B(x, 
\varepsilon ))}\right)^{\alpha \beta} 
{d \varepsilon \over \varepsilon})^{1/\beta}, &&(1.4)\cr}$$
and 
$$\gamma_{\alpha, \beta} (T, d) = \inf\limits_\mu \gamma_{\alpha, \beta} 
(T, d , \mu)$$
the infimum being taken on all probability measures.  Where $\beta = 1$, 
we write $\gamma_\alpha$ rather than $\gamma_{\alpha, 1}$.  When $T$ is a 
subset of a Hilbert space, $d$ will be the distance induced by the norm, and 
we will write simply $\gamma_{\alpha, \beta} (T)$.  Thus $\gamma_{1/2, 1}
(T) = 
\gamma_{1/2}(T)$.   The motivation for the  
introduction of these functionals is not a desire of empty generality, 
but the existence of 
concrete situations where these functionals are easy to
manipulate.  It is shown in [T7] how to compute $\gamma_{\alpha,\beta}
(\cl E)$ when $\cl E$ is an ellipsoid in Hilbert space, and it is shown 
that this computation is at the root of deep matching theorems of Ajtai, 
Komlos, Tusnady and Leighton and Shor on random samples in $[0,1]^2$.  This 
computation is done in [T7] by an explicit construction.  We will show that 
what actually only matters is the fact that the ellipsoid is $2$-convex.  
Equally irrelevant is the fact that we try to cover the ellipsoid with 
balls of a Hilbert space.  Only the covering numbers are relevant.  We 
recall that for two convex sets $B,U$ in a vector space, $N(B,U)$ denotes 
the minimum number of translates of $U$ needed to cover $B$.

\proclaim{Theorem 1.1}  Consider balanced convex sets $B,U$ in a vector 
space, and denote by
$\Vert\cdot\Vert_B$, $\Vert \cdot\Vert_U$ their gauges.  Assume that 
$$\lno{\Vert\cdot\Vert_B~~\text{is}~~2-\text{convex}.&&(1.5)\cr}$$

Then, for any $\alpha>0$, 
$$\lno{\gamma_{\alpha,2}(B,\Vert\cdot\Vert_U)\le
K\sup\limits_{\varepsilon > 0}\varepsilon (\log N (B, \varepsilon 
U))^\alpha, &&(1.6)\cr}$$ 
where $K$ depends only on
 the constant implicit in (1.5) and on $\alpha$.
\endproclaim

It is also easier in that case to get convinced that the left-hand side of
(1.6) should a priori be of bigger order than the right-hand side.  In the 
case of the ellipsoids of [T7], the computation of the right-hand side 
of (1.7) is a standard exercise (by volume estimates).  We should also 
mention 
that there is a converse to Theorem 1.1 when one replaces (1.5) by the 
condition that $\Vert\cdot\Vert_B$ is $2$-smooth.  We will not prove it.

There are good reasons to believe that the scheme of proof common 
to Theorem 1.1 and to the majorizing measure theorem is
conceptually the correct approach.  But our attempts to use this approach 
on the most interesting classes of functions considered in [T7] have failed 
due to impassable technical (and combinatorial) difficulties.  Nonetheless 
in Section 7 we will demonstrate how to use that scheme to study a non 
trivial situation (that could also be handled by the methods of [T7]) with 
the belief that our approach contains some of the essential ingredients 
needed for the final solution.

The rest of the paper is devoted to Bernoulli processes.  (The link with 
the previous material being that the crucial ingredients in each of the 
results we will present are closely related to the basic scheme of Section 
2.)

Consider an independent sequence $(\varepsilon_i)_{i\le M}$ of Bernoulli
random variables, i.e. $P(\varepsilon_i=1)=P(\varepsilon_i=-1)=1/2$.  For 
a subset $T$ of $\Bbb R^M$, we set  
$$\lno{b(T)=\sup\{E\sup\limits_{t\in S}\sum_{i\le M}\varepsilon_it_i;~~S
\subset T,~~S~~\text{finite}\}.&&(1.7)\cr}$$
(The ``Bernoulli process'' is the
collection of random variables $X_t=\sum\limits_{i\le M}\varepsilon_it_i$). 
The importance of this quantity is that sequences of ones and minus ones 
are one of the oldest and most fundamental structures of Probability.  How 
to describe $b(T)$ the way (1.1) describes $G(T)$?  A first simple 
observation is that, by a simple comparison result, we have 
$$\lno{b(T)\le KG(T).&&(1.8)\cr}$$
There, as well as in the rest of the paper, $K$ denotes a
universal constant, not necessarily the same at each occurrence.  (When the
constant depends only on, say, $\alpha,\beta$, we will write $K(\alpha,
\beta)$, etc.).

Another simple observation is that, using the bound $\vert\sum\limits_{i\le
M}\varepsilon_it_i\vert\le \sum\limits_{i\le M}\vert t_i\vert$, we have 
$$\lno{b(T)\le \sup\limits_{t\in T}\Vert t\Vert_1&&(1.9)\cr}$$
where $\Vert t\Vert_1=\sum\limits_{i\le M}\vert t_i\vert$.  Thus (1.8) and
(1.9) represent two different ways to control $b(T)$.  One can interpolate
between these two bounds, i.e.  $$\lno{b(T)\le K\inf\{u>0;T\subset
U+uB_1;~~~G(U)\le u\}&&(1.10)\cr}$$ where $B_1=\{t\in \Bbb
R^M;\sum\limits_{i\le M}\vert t_i\vert\le 1\}$.

Thus, it is natural to ask whether the right-hand side of (1.10) is actually
of the same order as the left hand side.

\subhead The Bernoulli Problem\endsubhead  Is it true that there exists a 
universal constant $K$ such that, given a subset $T$ of $\Bbb R^M$, one can 
find a subset $U$ of $\Bbb R^M$ such that $T\subset U+Kb(T)B_1$, where 
$\gamma_{1/2}(U)\le Kb(T)$?

(We have replaced in this statement $G(U)$ by the equivalent quantity
$\gamma_{1/2}(U)$).  The difficulty is of course that the decomposition is
neither unique nor canonical.  The Bernoulli Problem has been the main
motivation behind the papers [T4] [T5] that study related questions, some 
of them being discussed in Section 4, where the Bernoulli Problem will be
commented in detail.

The control on $\gamma_{1/2}(U)$ involves the structure of $U$ for the
$\ell_2$ norm, and the difficulty of the Bernoulli problem is that we have 
to separate a structure involving the $\ell_2$ norm with one involving the
$\ell_1$ norm.  In the special situation where on has a control on the
$\ell_{\infty}$ norm, the $\ell_1$ part disappears, making the problem 
easier.  In Section 4, we will explain why the correct formulation of this 
phenomenon is as follows.

\proclaim{Theorem 1.2} $\gamma_{1/2}(T)\le K(b(T)+\gamma_1(T,\Vert\cdot
\Vert_{\infty}))$.   \endproclaim

A basic tool in the proof of Theorem 1.2 is an extension of the construction
scheme of Section 2 to a two parameter situation.  This construction 
incidentally allows to recover some of the most technical results of [T5] 
with a significantly simpler proof.  This tool is presented in Section 5.

It turns out, (somewhat unexpectedly) that Theorem 1.2 is the starting point 
for an apparently new comparison principle between gaussian and Bernoulli 
averages.

\proclaim {Theorem 1.3}  Consider vectors $(x_i)_{i \le M}$ is a Banach $X$ 
of dimension $n$.  Then there is a subset $I$ of $\{1, \cdots , M\}$ 
such that 
$$\card I \le K n \log n \log {E\Vert \sum_{i \le M} g_i x_i \Vert 
\over E \Vert \sum_{i \le M} \varepsilon_i x_i \Vert} \le K n \log n 
\log\log n$$
such that 
$$E \Vert \sum_{i \not\in I , i \le M} g_i x_i \Vert \le K E \Vert  
\sum_{i \le M} \varepsilon_i x_i\Vert.$$
\endproclaim

While we have been unable to solve the Bernoulli problem, we did succeed in
proving the following weaker result.

\proclaim{Theorem 1.4}  Consider $p>1$.  Then there is a constant $K(p)$, 
that depends on $p$ only,
such that for each subset $T$ of $\Bbb R^M$, we can find a subset $U$ with
$T\subset U+K(p)b(T)B_p$, where $\gamma_{1/2}(U)\le K(p)b(T)$ and 
$$B_p=\{t\in\Bbb R^M;\sum_{i\le M}\vert t_i\vert^p\le 1\}.$$
\endproclaim

We do not know whether in Theorem 1.4 it is possible to replace $B_p$ by the
weak-$\ell_1$, ball
$$B_{1,\infty}=\{t\in\Bbb R^{\Bbb N};\quad\sup\limits_{u>0}u\card \{i;\vert 
t_i\vert\ge u\}\le 1\}$$
(a question that is apparently easier that the full Bernoulli problem).  
Apparently, all the available techniques are hopelessly inadequate to 
approach this question.  Thus it seems of interest to state a particularly 
attractive special case (that was pointed out to me by 
S. Montgomery-Smith).  Consider a finite group $G$, a function $f$ on $G$ 
such that $\Vert f \Vert_\infty \le 1$, and 
$$E \sup\limits_{t \in G} \vert \sum_{x \in G} \varepsilon_x f (t x) \vert 
\le 1.$$
There $(\varepsilon_x)_{x \in G}$ denotes an independent Bernoulli 
sequence, and $tx$ denotes the product in $G$.  For $i \ge 0$, consider 
the distance $d_i$ on $G$ given by
$$d^2_i (s, t) = \sum_{ x \in G} \min (2^{-4i}, (f (s x) - f(tx))^2)$$
The question is to decide whether, for some universal constant $K$, we have 
$$\sum_{i \ge 0} 2^{-i} \sqrt {\log N(G, d_i, 2^{-i})} \le K.$$
(The relevance of this inequality is not obvious and will be explained in 
Section 4).

The importance of Theorem 1.4 is that, while this result does not describe 
the structure of bounded Bernoulli processes, it is sufficient to describe 
the Radmacher cotype 2 constant of operators from $C(L)$ (where $L$ is a 
compact space).  Let us recall that the Rademacher cotype $2$ constant 
$C^r_2(V)$ of an operator $V$ from $C(L)$ to a Banach space $X$ is the 
infimum of the numbers $A$ such, for each $M$, and each sequence $(f_i)_{i
\le M}$ of continuous functions on $L$, we have  
$$\lno{\left(\sum_{i\le M}\Vert V(f_i)\Vert^2\right)^{1/2}\le AE\Vert
\sum_{i\le M}\varepsilon_if_i\Vert_{\infty}.&&(1.11)\cr}$$ 
Similarly, one defines the Gaussian
cotype $2$ constant $C^g_r(V)$ of $V$ as the infimum of the numbers $A$ 
such that for each $M$ and each sequence of continuous functions $(f_i)_{i
\le M}$ on $L$ we have 
$$\lno{(\sum_{i\le M}\Vert V(f_i)\Vert^2)^{1/2}\le
AE\Vert\sum_{i\le M}g_if_i\Vert&&(1.12)\cr}$$
Recall also that the $(2,1)$-summing norm $\Vert V\Vert_{2,1}$ of $V$ is 
defined as the infimum of the numbers $A$ such that for each $M$ and each 
sequence $(f_i)_{i\le M}$ of $C(L)$,  
$$\lno{(\sum_{i\le M}\Vert V(f_i)\Vert^2)^{1/2}\le A\Vert\sum_{i\le
M}\vert f_i\vert~\Vert_{\infty}.&&(1.13)\cr}$$

\proclaim{Theorem 1.5}  For some universal constant $K$, we have 
$$\lno{K^{-1}(C^g_2(V)+\Vert V\Vert_{2,1})\le C^r_2(V)\le K(C^g_2(V)+\Vert
V\Vert_{2,1}).&&(1.14)\cr}$$
\endproclaim

To help the reader to appreciate this result, we should state the 
following corollary.

\proclaim{Corollary 1.6}  Consider an operator $V$ from $\ell^N_{\infty}$ 
to a Banach space $X$.  Then 
$$\lno{C^r_2(V)\le K(LLN)^{1/2}\Vert V\Vert_{2,1},&&(1.15)\cr}$$
where $LLN=\max(1,\log(\log N))$.
\endproclaim

This result is optimal, and represents an improvement of order $(LLN)^{1/2}$
over the earlier result of [MS-T].

As a last result, we will mention the following corollary of Theorem 1.3.

\proclaim {Corollary 1.7}  In a Banach space $X$, of dimension $n$, there 
exists vectors $(x_i)_{i \le N}$ where 
$$N \le K n \log n LLn$$
and
$$C^r_2(X) E \Vert \sum_{i \le N} \varepsilon_i x_i \Vert < K 
( \sum_{i \le N} \Vert x_i \Vert^2)^{1/2}.$$
\endproclaim

In other words, the rademacher cotype $2$ constant of $X$ can be 
``computed with $N$ vectors'' within a universal constant.

Having stated the most easily understood results, we now describe in 
detail
the organization of the paper.  Section 2 gives the basic partitioning 
scheme and explains its relationship with the construction of majorizing 
measures.  As a first application we prove the Majorizing Measure Theorem 
for Gaussian processes.  In Section 3, the scheme of Section 2 is used to 
prove (a generalization of) Theorem 1.1.  The study of Bernoulli processes 
starts in Section 4.  
There the main problem is reformulated in several ways, that are shown to 
be equivalent.  The structure of Bernoulli processes does not depend (in 
contrast with the Gaussian case) of one single distance.  We recall the 
method (introduced in [T5]) of measuring the size of a set, provided with 
a family of distance-like functionals, with respect to the existence of 
majorizing measures.  In the rest of Section 4, we construct two families 
of functionals such that the resulting measure of size coincide exactly 
with the measure of size occurring in the Bernoulli problem (resp. a 
weaker form of the Bernoulli problem).  The motivation for these rather 
technically difficult results is that they demonstrate that the type of 
decomposition occurring in the Bernoulli problem is not intractable.  
However, these results are (mostly) not used in the sequel and should be 
omitted at first reading.  In Section 5, we establish the basic
partitioning scheme in the case of a family of distances.  This method 
extends, and significantly simplifies, crucial sections of the paper [T5].  
This scheme is then applied to the proof of Theorem 1.2, which in turn 
is applied to the proof of Theorem 1.3.  In Section 6, we 
prove Theorem 1.4.  Unfortunately, it is apparently impossible to use the 
scheme of Section 5, so we develop a specialized method.  We also prove 
Corollary 1.6.  Section 7 depends 
only on the material of Section 2.  Its purpose is to demonstrate how to 
use the basic partitioning scheme to study a concrete class of functions.   
\bigskip

\subhead 2.~~~A general partitioning construction\endsubhead

We consider a metric space $(T,d)$.  We denote by $B(x,a)$ the ball 
centered at $x$ of radius $a$.  For a subset $S$ of $T$, the diameter 
$\Delta(S)$ is defined as
$$\Delta(S)=\sup\{d(x,y);x,y\in S\}.$$
We assume that $T$ is of finite diameter.  Assume that, for 
$k\in\Bbb Z$, 
we are given a map $\varphi_k:T\to \Bbb R^+$.  We assume that 
$$\lno{&\forall\, x\in T,~~\forall\,k\in\Bbb Z,~~\varphi_k(x)\le
\varphi_{k+1}(x),&(2.1)\cr &A=\sup\{\varphi_k(x);~~x\in T,~~k\in\Bbb
Z\}<\infty.&(2.2)\cr}$$

Consider a function $\theta:\Bbb N\to \Bbb R^+$, and assume that 
$$\lno{\lim\limits_{n\to\infty}\theta(n)=\infty.&&(2.3)\cr}$$

Assume that, for certain numbers $r\ge 4$, $\beta>0$, the following 
holds:

\n (2.4)\quad Given any point $x$ of $T$, any $k\in\Bbb Z$, any $n\ge 1$, 
and any points $y_1,\dots, y_n$ of $B(x,r^{-k})$ such that 
$$\forall~i,j\le n,~~i\not= j\Rightarrow d(y_i,y_j)\ge r^{-k-1}$$
we have 
$$\max\limits_{j\le n}\varphi_{k+2}(y_j)\ge \varphi_k(x)+r^{-\beta k}
\theta(n).$$
The reader should observe that on the left we have $\varphi_{k+2}$ 
rather 
than $\varphi_{k+1}$.  This is the crucial point of the condition.

We denote by $k_0$ the largest integer such that the diameter of $T$ is 
$\le r^{-\beta k_0}$.

\proclaim{Theorem 2.1}  Under conditions (2.1) to (2.4), one can find an
increasing sequence of finite partitions $(\cl C_k)_{k\ge k_0}$ of $T$, 
and for each atom $C$ of $\cl C_k$ one can find an index $\ell_k (C)\ge 1$, 
such that the following properties hold.

\n (2.5)\quad Each set $C$ of $\cl C_k$ has diameter $\le 2\cdot r^{-k}$.

\n (2.6)\quad For any $k\ge k_0+1$, and any two sets $C,D$ of $\cl C_{k+1}$
that are included in the same element of $\cl C_k$, we have $\ell_{k+1}(C)
\not= \ell_{k+1}(D)$.

\n (2.7)\quad If, for $x\in T$, we denote by $C_k(x)$ the unique element of
$\cl C_k$ that contains $x$, we have 
$$\forall~x\in T,\quad \sum_{k\ge k_0}r^{-\beta k}\theta(\ell_{k+1}(C_{k+1}
(x)))\le 4A.$$
\endproclaim

\subhead Comment\endsubhead (2.7) is of course a ``smallness condition'' 
on the sequence $(\cl C_k)_{k\ge k_0}$.

\demo{Proof}  The construction goes by induction over $k$.  For each set
$C\in\cl C_k$, we also construct a distinguished point $z_k(C)\in C$, such 
that
$$\lno{\forall~y\in C,\quad d(y,z_k(C))\le r^{-k}.&&(2.8)\cr}$$
Observe that this condition implies (2.5).
\enddemo

At the first stage, we set $\cl C_{k_0}=\{T\},\ell_{k_0}(T) = 1$, and we 
choose $z_{k_0}(T)$ such that
$$\varphi_{k_0+2}(z_{k_0}(T))\le 2^{-1}A+\inf\{\varphi_{k_0+2}(y);y\in
T\}.$$

We assume now that the partition $\cl C_k$ has been constructed, as well as
the points $z_k(C)$ for $C\in\cl C_k$.  To construct $\cl C_{k+1}$, it 
suffices to show how to partition any given element $C$ of $\cl C_k$.  This 
will be done in turn by an inductive argument.  The index of any piece will 
simply be the rank at which it is constructed.

First, we choose $y_1\in C$ such that 
$$\varphi_{k+2}(y_1)\le 2^{k_0-k-1}A+\inf \{\varphi_{k+2}(y);y\in C\}.$$
We then pick inductively $y_2,\dots, y_{\ell}$ such that 
$$y_{\ell}\in C\bs \bigcup\limits_{i<\ell}B(y_i,r^{-k-1})$$
and 
$$\lno{\varphi_{k+2}(y_{\ell})\le 2^{k_0-k-1}A+\inf \{\varphi_{k+2}(y);y\in 
C\bs \bigcup\limits_{i<\ell}B(y_i,r^{-k-1})\}.&&(2.9)\cr}$$

The construction continues as long as possible.  It eventually stops 
according to (2.2), (2.3), (2.4).  (This point will be detailed later.)  
We set 
$$D_{\ell}=(C\cap
B(y_{\ell},r^{-k-1}))\bs\bigcup\limits_{i<\ell}B(y_i,r^{-k-1}).$$
Let us observe immediately the following crucial fact, that follows from 
(2.9)
$$\lno{\forall~y\in D_{\ell},\quad \varphi_{k+2}(y_{\ell})\le
2^{k_0-k-1}A+\varphi_{k+2}(y).&&(2.10)\cr}$$

The sets $D_{\ell}$ form the partition of $C$ that we look for.  We set 
$$\lno{&z_{k+1}(D_{\ell})=y_{\ell}&(2.11)\cr
&\ell_{k+1}(D_{\ell})=\ell.&(2.12)\cr}$$
Thus (2.10) can be rewritten as 
$$\lno{\forall~y\in D_{\ell},\quad \varphi_{k+2}(z_{k+1}(D_{\ell}))\le
2^{k_0-k-1}A+\varphi_{k+2}(y).&&(2.13)\cr}$$

This completes the construction.  It is obvious that (2.8) (hence (2.5)) 
holds, and we proceed to prove (2.7).

We observe that $d(z_k(C),y_i)\le r^{-k}$, and that for $i<j$, we have
$d(y_i,y_j)\ge r^{-k-1}$.  Thus, by (2.4), for each $\ell$ we have 
$$\lno{\max\limits_{i\le \ell}\varphi_{k+2}(y_i)\ge \varphi_k(z_k(C))+r^{-
\beta k}\theta(\ell).&&(2.14)\cr}$$

\n On the other hand, by (2.10), for $i\le \ell$, and since $y_{\ell}\in 
D_i$, we have 
$$\varphi_{k+2}(y_i)\le 2^{k_0-k-1}A+\varphi_{k+2}(y_{\ell})$$
so that, combining with (2.14) we get 
$$\lno{\varphi_{k+2}(y_{\ell})+2^{k_0-k-1}A\ge \varphi_k(z_k(C))+r^{-\beta
k}\theta(\ell).&&(2.15)\cr}$$

Consider now any $x\in T$, and apply (2.15) to $C=C_k(x)$, $\ell=\ell_{k+1}
(C_{k+1} (x))$, $D_{\ell}=C_{k+1}(x)$, $y_{\ell}=z_{k+1}(D_{\ell})=z_{k+1}
(C_{k+1}(x))$, so that 
$$\lno{2^{k_0-k-1}A+\varphi_{k+2}(z_{k+1}(C_{k+1}(x)))\ge
\varphi_k(z_k(C_k(x)))+r^{-k\beta}\theta(\ell_{k+1}(C_{k+1}(x))).&&(2.16)
\cr}$$
Since 
$$z_{k+2}(C_{k+2}(x))\in C_{k+2}(x)\subset C_{k+1}(x)=D_{\ell},$$
by (2.13) we have 
$$\varphi_{k+2}(z_{k+1}(C_{k+1}(x)))\le 2^{k_0-k-1}A+\varphi_{k+2}(z_{k+2}
(C_{k+2}(x))).$$
Combining with (2.16) we get 
$$\lno{2^{k_0-k}A+\varphi_{k+2}(z_{k+2}(C_{k+2}(x)))\ge
\varphi_k(z_k(C_k(x)))+r^{-k\beta}\theta(\ell_{k+1}(C_{k+1}(x))).&&(2.17)
\cr}$$
If we sum these relations for $k_0\le k\le m$, we obtain
$$\eqalign{\sum_{k_0\le k\le m}r^{-k\beta}\theta(\ell_{k+1}(C_{k+1}(x)))
&\le 2A+\varphi_{m+1}(z_{m+1}(C_{m+1}(x))))\cr
&\quad +\varphi_{m+2}(z_{m+2}(C_{m+2}(x)))\le 4A\cr}$$
by (2.2).  This completes the proof.\hfill\bx

The following result clarifies the relationship between the situation of
Theorem 2.1 and majorizing measures. 

\proclaim{Theorem 2.2}  Assume that we have an increasing sequence of
finite partitions $(\cl C_k)_{k\ge k_0}$ on $T$, and assume that for each 
atom $C$ of $\cl C_k$, we have an index $\ell_k(C)\ge 1$, that satisfies 
(2.6).  Assume moreover that for some numbers $\alpha>0$, $M$, we have 
$$\lno{\forall~x\in T,\quad \sum_{k\ge k_0}r^{-\beta k}(\log
\ell_{k+1}(C_{k+1}(x)))^{\alpha}\le M,&&(2.18)\cr}$$ 

\n where $C_k(x)$ denotes again the
unique element of $\cl C_k$ that contains $x$.

Then, we can find a probability measure $\mu$ on $T$ such that 
$$\lno{\forall~x\in T,\quad \sum_{k\ge k_0}r^{-\beta
k}\left(\log{1\over\mu(C_k(x))}\right)^{\alpha}\le K(\alpha,\beta, r)(M+
r^{-\beta k_0})&&(2.19)\cr}$$
where $K(\alpha,\beta,r)$ depends on $\alpha,\beta,r$ only.
\endproclaim

\demo{Proof}  We first observe the elementary fact that 
$$\lno{\sum_{\ell\ge 1}{1\over (\ell+1)^2}\le 1.&&(2.20)\cr}$$
\enddemo

We set $w_{k_0}(T)={1\over 2}$, and, for $k>k_0$ we define inductively a
number $w_k(C)$ for $C\in \cl C_k$ by 
$$\lno{w_k(C)={1\over 2(\ell_k(C)+1)^2}w_{k-1}(C')&&(2.21)\cr}$$
where $C'$ is the element of $\cl C_{k-1}$ that contains $C$.  Using (2.20)
and (2.6) we see inductively that 
$$\sum_{C\in\cl C_k}w_k(C)\le 2^{k_0-k-1}.$$
Thus, we can find a probability measure $\mu$ on $T$ that gives mass $\ge
w_k(C)$ to an arbitrary point $z_k(C)$, of $C$, for $k\ge k_0$, $C\in
\cl C_k$.  Thus we have 
$$\sum_{k>k_0}r^{-\beta k}\left(\log{1\over\mu(C_k(x))}\right)^{\alpha}\le
H=:\sum_{k>k_0}r^{-\beta k}\left(\log{1\over w_k(C_k(x))}\right)^{\alpha}.$$
Using (2.21), we get 
$$\lno{H\le \sum_{k>k_0}r^{-\beta k}\left[\log{1\over
w_{k-1}(C_{k-1}(x))}+\log(2(\ell(C_k(x))+1)^2)\right]^{\alpha}.&&(2.22)
\cr}$$
Now, we observe that, if we set $\delta=(1+r^{\beta})/2$, for all
$x,y>0$, we have, since $\delta>1$  $$(x+y)^{\alpha}\le \delta
x^{\alpha}+K(\alpha,\beta,r)y^{\alpha}.$$

Thus, from (2.22), we get 
$$H\le \delta\sum_{k>k_0}r^{-\beta k}\left(\log{1\over
w_{k-1}(C_{k-1}(x))}\right)^{\alpha}+K(\alpha,\beta,r)\sum_{k>k_0}
r^{-\beta
k}(\log 2(\ell(C_k(x))+1)^2)^{\alpha}.$$
Since the first summation is at most 
$$\delta r^{-\beta k_0}(\log 2)^{\alpha}+\delta r^{-\beta}H$$
and $\delta r^{-\beta}\le 1-(1-r^{-\beta})/2$, the result follows easily.
\hfill\bx

\remark {Remark}  When $\alpha \le 1$, using that fact that $(x + y)^\alpha 
\le x^\alpha + y^\alpha$, we see that we can take $K(\alpha, \beta, r)= K$, 
independent of $\alpha, \beta , r$.
\endremark

As a first application, we prove the left-hand side of (1.3).  For this, 
we use Theorem 2.1 with $\theta(n)=K^{-1}\sqrt{\log n}$ (for a large enough 
$K$) and
$$\varphi_k(x)=G(T)-G(B(x,r^{-k})).$$
The fact that (2.4) holds for $r$ sufficiently large is proved in [T4].  
The left-hand side of (1.3) then follows from Theorem 2.2.  It is of 
interest to compare the approach of Theorem 1.1 with that of [T4].  The 
main idea is identical but the argument is, so to say, reversed.  One gain 
in this approach is that we no longer need the analogue of Lemma of [T6]; 
this is fortunate, since this analogue would not hold for all the values of 
$\alpha,\beta$ of interest.

\bigskip

\subhead 3.~~Majorizing measures on sufficiently convex sets\endsubhead

In this section $T=B$ is the unit ball of a normed space $X$.  We assume 
that the norm of $X$ has a modulus of convexity with a $p^{\text{th}}$ 
power estimate.  More precisely we assume that for some number $p$ $(p\ge 2
)$ and some number $\gamma>0$, we have  
$$\lno{\inf \left\{1-{\Vert x+y\Vert\over 2};\Vert
x\Vert=\Vert y\Vert=1;~~\Vert x-y\Vert\le \varepsilon \right\}\ge \gamma
\varepsilon ^p,&&(3.1)\cr}$$
for all $0<\varepsilon\le 2$. The choice of parameters $(p = \beta = 2)$ 
made in the statement of Theorem 1.1 is important for applications, but 
it is 
instructive (and require no further effort) to perform the proof in a more 
general setting.

Let us first note the following simple fact 

\proclaim{Lemma 3.12}  If $\Vert x\Vert$, $\Vert u\Vert\le t$ and $\Vert
x+y\Vert\ge 2u$, then
$$\gamma \Vert x-y\Vert^p\le t^{p-1}(t-u).$$
\endproclaim

\demo{Proof}  In [L-T], p. 60 it is shown that in (3.1) one can replace the
condition $\Vert x\Vert=\Vert y\Vert=1$ by $\Vert x\Vert$, $\Vert y\Vert\le
1$.  If we use (3.1) for $x/t$, $y/t$, we then obtain 
$$1-{\Vert x+y\Vert\over 2t}\ge \gamma{\Vert x-y\Vert^p\over t^p}$$
from which the result follows.\hfill\bx
\enddemo

The proof of Theorem 3.1 consists of two main steps.  In the first, we will
apply Theorem 2.1, and in the second we will apply Theorem 2.2.

To apply Theorem 2.1, we define the functions $\varphi_k(x)$ for $x\in B$ 
as follows.  We set $\varphi_k(x)=0$ if $x\in 2r^{-k}U$.  Otherwise, we set
$$\varphi_k(x)=\sup\{t>0;tB\cap (x+2r^{-k}U)=\emptyset\}.$$
It is obvious that (2.1) and (2.2) hold, with $A=1$.  

For $n\ge 2$, we set 
$$\varepsilon(n)=\sup\{\varepsilon>0;\exists~y_1,\dots, y_n\in B;
\forall~\ell,~\ell',~1\le \ell<\ell'\le n,y_{\ell}\not\in y_{\ell'}+
\varepsilon U\}.$$
We observe that $\varepsilon(n)\le 2$.  We observe the following simple
relations  
$$\lno{\varepsilon<{\varepsilon(n)\over 2}\Rightarrow &N\left(B,\varepsilon 
U\right)\ge n &(3.3)\cr
&\varepsilon'>\varepsilon(n)\Rightarrow N(B,\varepsilon'U)\le n.&(3.4)\cr}$$

\proclaim{Lemma 3.3}  Assume $r\ge 8$.  Consider $x\in B$, $n\ge 2$, and 
points $y_1,\dots, y_n\in B(x,r^{-k})$, such that 
$$i<j\le n\Rightarrow y_j\not\in y_i+r^{-k-1}B.$$
Then we have 
$$\sup\limits_{\ell\le n}\varphi_{k+2}(y_{\ell})\ge
\varphi_k(x)+{\gamma r^{-kp}\over (2r\varepsilon(n))^p}.$$
\endproclaim

\demo{Proof}  By definition of $\varphi_{k+2}$, for $t>\sup\limits_{i\le
n}\varphi_{k+2}(y_i)$ and any $i\le n$, we can find a point 
$$z_i\in tB\cap (y_i+2r^{-k-2}U).$$
We note that 
$$z_i\in y_i+2r^{-k-2}U\subset x+(2r^{-k-2}+r^{-k})U \subset
x+2r^{-k}U.$$
\enddemo

Thus, by convexity of $U$, we have 
$$\forall~i,j\le n,\quad {z_i+z_j\over 2}\in x+2r^{-k}U$$
and thus, by definition of $\varphi_k$, we have 
$$\Vert {z_i+z_j\over 2}\Vert \ge \varphi_k(x).$$

Since $z_i\in tB$ for all $i\le n$, it follows from Lemma 3.1 that 
$$\gamma\Vert z_i-z_j\Vert^p\le t^{p-1}(t-\varphi_k(x)).$$

It thus follows that the points $z_i$, for $i\le n$, belong to the ball
$z_1+RB$, where $R^p=(t-\varphi_k(x))/\gamma$.  

Now, since $z_i\in y_i+2r^{-k-2}U$, since $y_i\not\in r^{-k-1}U+y_j$ for
$j\not= i$, and since $r\ge 8$, we have 
$$ z_i-z_j\not\in {r^{-k-1}\over 2}U.$$
By definition of $\varepsilon(n)$, we thus have 
$${r^{-k-1}\over 2R}\le \varepsilon(n)$$
which means 
$$t^{p-1}(t-\varphi_k(x))\ge {\gamma\over (2r)^p}{r^{-kp}\over \varepsilon
(n)^p}.$$
Since $t>\sup\limits_{\ell\le n}\varphi_{k+2}(y_{\ell})$ is arbitrary, and 
since $A=1$, this completes the proof.\hfill\bx

We now can apply Theorem 2.1 with $\theta(1)=0$,
$\theta(n)=\gamma/(2r\varepsilon(n))^p$ for $n\ge 2$.  We observe that we 
have $k_0=0$.

\proclaim{Corollary 3.4}  We can find an increasing sequence of finite
partitions $(\cl C_k)_{k\ge 0}$ of $B$, and indexes $\ell_k(C)$ for $C\in 
\cl C_k$, that satisfy conditions (2.5) and (2.6) of Theorem 2.1, and such 
that 
$$\lno{\forall x\in B,\quad \sum_{k\ge 0}{r^{-kp}\over \varepsilon(
\ell_{k+1}(C_{k+1}(x)))^p}\le {2^{p+2}r^p\over \gamma}&&(3.5)\cr}$$
where we make the convention that, when $\ell_{k+1}(C_{k+1}(x))=1$, the
corresponding term of the series is zero. \endproclaim

The next main step in the proof of Theorem 3.1 is to interpret condition 
(3.5) when we suitably control the covering numbers $N(B,\varepsilon U)$.

\proclaim{Proposition 3.5}  Consider the number $\beta'>0$ such that 
$$\lno{{1\over\beta}={1\over\beta'}+{1\over p}&&(3.6)\cr}$$
and assume that 
$$\lno{M^{\beta'}=\sum_{k\ge 0}r^{-k\beta'}(\log
N(B,r^{-k}U))^{\alpha\beta'}<\infty.&&(3.7)\cr}$$ 

Then the partition of Corollary 3.4 satisfies 
$$\lno{\sum_{k\ge
0}r^{-k\beta}(\log\ell_{k+1}(C_{k+1}(x)))^{\alpha\beta}\le K(\alpha,\beta, 
r) {M^{\beta}\over \gamma^{\beta/p}}.&&(3.8)\cr}$$
\endproclaim

\subhead Comment\endsubhead In (3.6) we certainly allow the case
$\beta'=\infty$ $(\beta=p)$, in which case (3.7) has to be interpreted as 
$$M=\sup\limits_{k\ge 0}r^{-k}(\log N(B,r^{-k}U))^{\alpha}<\infty.$$
The necessary modifications to the proof in that case are left to the 
reader.

\demo{Proof}  {\bf Step 1.}\quad We fix $x$.  We observe that in (3.8) the
contribution of the terms for which $\ell_{k+1}(C_{k+1}(x))=1$ is zero.  
For $m\ge 0$, we set
$$I(m)=\{k\ge 0;2^m\le\log_2\ell_{k+1}(C_{k+1}(x))<2^{m+1}\}.$$
When $I(m)$ is not empty, we denote
by $i(m)$ its smallest element, and we make the convention that when $I(m)$ 
is empty the corresponding term does not appear.  We have
$$\eqalign{\sum_{k\in I(m)}r^{-k\beta}(\log_2\ell_{k+1}(C_{k+1}(x)))^{\alpha
\beta}&\le \sum_{k\in I(m)}r^{-k\beta}2^{(m+1)\alpha\beta}\cr
&\le \sum_{k\ge i(m)}r^{-k\beta}2^{(m+1)\alpha\beta}\cr
&\le K(\alpha,\beta)r^{-i(m)\beta}2^{m\alpha\beta}.\cr}$$
Thus 
$$\lno{\sum_{k\ge 0}r^{-k\beta}(\log\ell_{k+1}(C_{k+1}(x)))^{\alpha\beta}
\le K(\alpha,\beta)\sum_{m\ge 0}r^{-i(m)\beta}2^{m\alpha\beta}.&&(3.9)\cr}$$

\n {\bf Step 2.}\quad We set $u=\beta/\beta'$, $v=\beta/p$, so that, by 
(3.6), $u+v=1$.  We observe the identity 
$$r^{-i(m)\beta}2^{m\alpha\beta}=(2^{m\alpha\beta'}
\varepsilon^{\beta'}_m)^u(r^{-pi(m)}\varepsilon^{-p}_m)^v$$
where $\varepsilon_m=\varepsilon(2^{2^m})$.
\enddemo

Thus, by H\"older's inequality, we have 
$$\lno{S\le S^u_1S^v_2,&&(3.10)\cr}$$
where 
$$S=\sum_{m\ge 0}r^{-i(m)\beta}2^{m\alpha\beta};~~S_1=\sum_{m\ge
0}2^{m\alpha\beta'}\varepsilon^{\beta'}_m;S_2=\sum_{m\ge
0}r^{-pi(m)}\varepsilon^{-p}_m.$$

\n{\bf Step 3.}\quad Since 
$$\ell_{i(m)+1}(C_{i(m)+1}(x))\ge 2^{2^m}$$
and since $\varepsilon(n)\ge \varepsilon(n')$ for $n\le n'$, we have 
$$\varepsilon(\ell_{i(m)+1}(C_{i(m)+1}(x)))^{-p}\ge \varepsilon^{-p}_m.$$
Since, from the definition, we have $i(m)\not= i(m')$ for $m\not= m'$, it
follows from (3.5) that
$$\lno{S_2\le {2^{p+2}r^p\over \gamma}. &&(3.11)\cr}$$

\n{\bf Step 4.}\quad Control of $S_1$.

For $i\ge 0$, consider the set 
$$J(i)=\{m\ge 0;2r^{-i-1}< \varepsilon_m\le 2r^{-i}\}.$$
When $J(i)$ is non-empty, we denote by $m(i)$ its largest element, and when
$J(i)=\emptyset$, we make the convention that the corresponding term does 
not appear.

We have 
$$\lno{\sum_{m\in J(i)}2^{m\alpha\beta'}\varepsilon_m^{\beta'}&\le
\sum_{\ell\le m(i)}2^{\ell\alpha\beta'}2^{\beta'}r^{-\beta'i}&(3.12)\cr
&\le K(\alpha,\beta)2^{m(i)\alpha\beta'}r^{-\beta'i}.\cr}$$

Since $\varepsilon_{m(i)}> 2r^{-i-1}$, by (3.2) we have 
$$N(B,r^{-i-1}U)\ge 2^{2^{m(i)}}$$  
since $\varepsilon_m=\varepsilon(2^{2^m})$.  Thus, using (3.7) we get  
$$\eqalign{\sum_{i\ge 0}2^{m(i)\alpha\beta'}r^{-\beta'i}&
\le \sum_{i\ge 0}r^{-\beta'i}(\log_2N(B,r^{-i-1}U))^{\alpha\beta'}\cr
&\le r^{\beta'}\sum_{i\ge 0}r^{-\beta'i}(\log_2N(B,r^{-i}U))^{\alpha
\beta'}\cr
&\le r^{\beta '} K(\alpha,\beta)M^{\beta'}.\cr}$$

Combining with (3.12), we get 
$$\lno{S_1\le r^{\beta'} K(\alpha,\beta)M^{\beta'}.&&(3.13)\cr}$$

The result then follows from (3.9) to (3.13).\hfill \bx

Theorem 3.1 is now a consequence of Theorem 2.2 and Proposition 3.5.
\bigskip

\subhead 4.~~Functionals on classes of functions and the Bernoulli
Problem\endsubhead

For $\tau\ge1$, let us consider the measure $\mu_{\tau}$ of density 
$a_{\tau}e^{-\vert t\vert^{\tau}}$ with respect to Lebesgue measure (where 
$a_{\tau}$ is a normalizing constant) and let us consider a family $(h_i)_{i
\le N}$ of independent random variables distributed like $\mu_{\tau}$.  For 
a subset $T$ of $\Bbb R^N$, we can consider the ``canonical process'' 
$(X_t)_{t\in T}$, where 
$$X_t=\sum_{i\le N}t_ih_i$$
and the quantity 
$$F_{\tau}(T)=E\sup\limits_{t\in T}X_t=:\sup\limits_{S\subset T,S~
\text{finite}}E\sup\limits_{t\in S}\sum_{i\le N}t_ih_i.$$

It is a remarkable fact that the quantity $F_{\tau}(T)$ can be characterized 
in terms of the geometry of $T$.  This is a generalization of the majorizing
measure theorem.  In the following we write $A\sim B$ to mean $A\le K(\tau)
B$, $B\le K(\tau)A$.  We denote by $\tau'$ the conjugate exponent of $\tau$.

\proclaim{Theorem 4.1}  [T6]

a)\quad If $1\le \tau\le 2$, we have
$$F_{\tau}(T)\sim\gamma_{1/2}(T,\Vert
\cdot\Vert_2)+\gamma_{1/\tau}(T,\Vert \cdot \Vert_{\tau'}).$$

b)\quad If $\tau\ge 2$, then 
$$F_{\tau}(T)\sim \min\{A;T\subset U+V;\gamma_{1/2}(T,\Vert \cdot\Vert_2)
\le A,\gamma_{1/\tau}(T,\Vert\cdot \Vert_{\tau'})\le A\}.$$
\endproclaim

The case $\tau=2$ is the majorizing measure theorem for Gaussian processes.   
Part b is strikingly similar with the Bernoulli Conjecture.  Actually this
conjecture can be thought of to be like the ``limiting case'' $\tau\to
\infty$ of Theorem 4.1.  (It should be noted that the constants implicit in 
the symbol $\sim$, as given by the arguments of [T6], go to infinity with 
$\tau$.)

Thus the majorizing measure theorem appears as one given element of a 
continuous family of theorems.  One may then wonder whether the Bernoulli 
conjecture is properly formulated, and why Gaussian processes should play 
a prominent part in that conjecture.  As it turns out, there is no reason 
to distinguish Gaussian process (other than their intrinsic importance).  
Actually, one could also choose to distinguish the canonical processes for 
$\tau=1$.  This is also natural, since, by comparison, we have $F_{\tau}(T)
\le K(\tau)F_{\tau'}(T)$ for $\tau\ge\tau'$, so that control of $F_1(T)$ is 
the strongest of this family of conditions.  As it turns out, these 
formulations are equivalent.

\proclaim{Theorem 4.2}  The following are equivalent:

a)\quad For any subset $T$ of $\Bbb R^N$, one can find $U\subset \Bbb R^N$ 
such that $T\subset U+Kb(T)B_1$, and $\gamma_{1/2}(U)\le Kb(T)$.

b)\quad For any subset $T$ of $\Bbb R^N$, one can find $U\subset \Bbb R^N$ 
such that $T\subset U+Kb(T)B_1$, and 
$$\gamma_{1/2}(U)\le Kb(T);\quad\gamma_1(U,\Vert\cdot\Vert_{\infty})\le 
Kb(T).$$

c)\quad For any subset $T$ of $\Bbb R^N$, one can find $U\subset \Bbb R^N$ 
such that $T\subset U+Kb(T)B_1$ and $\gamma_1(U,\Vert\cdot \Vert_{\infty})
\le Kb(T)$.
\endproclaim

\demo{Proof}  $a$) is the original formulation.  It is obvious that $b)
\Rightarrow a),c)$.
\enddemo

We prove that $c)\Rightarrow a)$.  Consider $U$ as given by $c)$ and 
consider the set 
$$V=U\cap (T+Kb(T)B_1).$$
Thus $b(V)\le Kb(T)$, and, since $V\subset U$, we have $\gamma_1(V,\Vert 
\cdot\Vert_{\infty})\le Kb(T)$.  Thus, by Theorem 1.2 (that we will prove 
in Section 5) we have $\gamma_{1/2}(V)\le Kb(T)$.  Now, since $T\subset U+
Kb(T)B_1$, it is easy to check that $T\subset V+Kb(T)B_1$, which proves 
a).

Thus, it remains to prove that $a\Rightarrow b$.  For this it suffices to 
show that for any subset $U$ of $\Bbb R^N$, we can write $U\subset W+aB_1$, 
where $a=K\gamma_{1/2}(U)$ and 
$$\gamma_{1/2}(W)\le Ka,\quad \gamma_1(W,\Vert \cdot\Vert_{\infty})\le K
a.$$

The idea to prove this goes back to [T1], and has subsequently been used
numerous times by this author.  We will not reproduce the argument,
since a stronger fact will be proved in Proposition 4.3 below.\hfill\bx

One of the difficulties of studying the Bernoulli problem is that the
corresponding measure 
$$\lno{\inf\{u;T\subset U+uB_1;~~\gamma_{1/2}(U)\le u\}&&(4.1)\cr}$$
of the size of a subset $T$ of $\Bbb R^N$ is very cumbersome to manipulate.  
The corresponding difficulty was solved in [T6] in the case of Theorem 4.1 
and of the canonical processes.  It is not obvious at all how to adopt 
these ideas to the case of Bernoulli processes.  We will nonetheless show 
that there exists indeed a seemingly more manageable functional that is 
exactly equivalent to the quantity (4.1).  This is certainly an encouraging 
fact.  

First we must recall a convenient notion (introduced in [T5]) that allows 
one to study the size of a set relatively to a family of functionals.  
Consider a set $T$, and assume that we are given a number $r\ge 4$ and that 
for $j\in\Bbb Z$ we have a function $\varphi_j$ on $T\times T$, $\varphi_j
\ge 0$, that satisfies $\varphi_j(s,t)=\varphi_j(t,s)$.  Typically 
$\varphi_j$ will be the square of a distance, so that we can assume 
$$\lno{\varphi_j(s,t)\le 4(\varphi_j(s,u)+\varphi_j(u,t))&&(4.2)\cr}$$
for all $s,t,u$ of $T$.  For a subset $U$ of $T$, we define 
$$D_j(U)=\sup\{\varphi_j(s,t);s,t\in U\}.$$

Consider $i\in \Bbb Z$ and an increasing sequence of finite partitions 
$(\cl C_j)_{j\ge i}$ of $T$.  For $x\in T$, $j\ge i$, we denote as usual by 
$C_j(x)$ the unique element of $\cl C_j$ that contains $x$.  For a 
probability measure $\mu$ on $T$ we consider the quantity
$$\sup\limits_{x\in
T}\sum_{j\ge i}r^{-j}\left(D_j(C_j(x))+\log{1\over\mu(C_j(x))}\right)$$
and we define the functional $\theta_i(T)$ as the infimum of the previous
quantity over all possible choices of $\mu$ and the sequence $(\cl C_j)_{j
\ge i}$.  (The reader should observe that, in contrast with the definition 
of [T5], we do {\it not} require that $\cl C_i=\{T\}$.  This is however 
only a minor technical point.)  The idea of these functionals,
as explained in [T5], Section 3 is that they are related to the usual 
notion of majorizing measures through a change of variable.  In particular 
we will use the following fact, that is proved in [T5] (and can also be 
deduced from Theorem 5.1 below).

\n (4.3)\quad When one uses the functionals
$$\varphi_j(s,t)=r^{2j}\Vert s-t\Vert^2_2$$
then 
$$\theta_i(T)\le K(r)(\gamma_{1/2}(T)+r^{-i}).$$

It will be notationwise more convenient to work now in the space of 
measurable functions on a measure space $(\Omega,\Sigma,\lambda)$.  We do 
not assume that $\lambda$ is a probability; indeed the most important case 
is $\Omega = \{1,\dots, N\}$, $\lambda$ being the counting measure.  The
more general formulation has also some intrinsic interest).  We denote by 
$B_1$ (resp. $B_2,B_{\infty}$) the unit ball of $L_1(\lambda)$ (resp.
$L_2(\lambda), L_{\infty}(\lambda)$).

On $L_2(\lambda)$, we consider the functions 
$$\lno{\varphi_j(f,g)=\int_{\Omega}\min(1,r^{2j}(f-g)^2)d\lambda.&
&(4.4)\cr}$$
An immediate observation is that 
$$\lno{\varphi_j(f,g)\le r^{2j}\Vert f-g\Vert^2_2.&&(4.5)\cr}$$ 

\proclaim{Proposition 4.3}  Consider a subset $T$ of $L_2(\lambda)$ and
$i\in\Bbb Z$.  Assume that $T\subset {r^{-i}}B_{\infty}/4$.  Then we have 
$$T\subset U+K\theta_i(T)B_1$$
where 
$$\lno{\gamma_{1/2}(U)\le
Kr(\theta_i(T)+r^{-i}+r^{-i+1}D_{i-1}(T));~~\gamma_1(U,\Vert
\cdot\Vert_{\infty})\le K(\theta_i(T)+r^{-i}).&&(4.6)\cr}$$
\endproclaim

\subhead Comments\endsubhead 1)\quad In the case where $\lambda$ the
counting measure, we have $B_1\subset B_2\subset B_{\infty}$.  For this 
choice of $\lambda$, the condition $T\subset r^{-i}B_{\infty}/4$ is not 
very restrictive in practice.  Indeed,  the Proposition will be used for 
values of $i$ such that $r^{-i}$ is of order $\theta_i(T)$.  For these 
values, and since $\diam_2 (U) \le K \gamma_{1/2}(U)$, (4.6) implies in 
any case that $T \subset K r^{-i} B_2$.

\n 2)\quad Clearly the value of $\theta_i(T)$ can only decrease when the
functionals $\varphi_j$ decrease.  Thus, by (4.3),  we see that $\theta_i
(T)\le K(\gamma_{1/2}(T)+r^{-i})$ so that Proposition 4.3 can be applied 
to complete the proof of $a\Rightarrow b$ in Theorem 4.2.

\demo{Proof}  {\bf Step 1.}\quad We start with a simple observation.  By
Markov's inequality, if $0\le a\le r^{-j}$, we have 
$$\lno{\lambda(\{\vert f\vert\ge a\})\le {1\over a^2}\int \min
(f^2,r^{-2j})d\lambda={r^{-2j}\over a^2}\int\min(r^{2j}f^2,1)d\lambda.
&&(4.7)\cr}$$

\n{\bf Step 2.}\quad Consider an increasing sequence $(\cl C_j)_{j\ge i}$ 
of finite partitions of $T$ and a probability measure $\mu$ on $T$ such 
that 
$$\lno{\forall~x\in T,\quad \sum_{j\ge i}r^{-j}\left(D_j(C_j(x))+\log{1
\over\mu(C_j(x))}\right)<2\theta_i(T).&&(4.8)\cr}$$
For $C\in \cl C_j$, $j\ge i$, we chose one element $y(C)\in C$.  We select 
one element $y(T)\in T$.  For $x\in T$, $j\ge i$, we set $\pi_j(x)=y(C_j(x)
)$.  We set $\pi_{i-1}(x)=y(T)$.

For each $\omega\in\Omega$, we define 
$$\lno{\ell(x,\omega)=\inf\{j\ge i-1;\mid\pi_j(x)(\omega)-\pi_{j+1}(x)(
\omega)\vert> r^{-j}\}. &&(4.9)\cr}$$
(When the set on the right is empty, we set  $\ell(x,\omega)=\infty)$.  
We set $\pi_{\infty}(x)=x$, and we set 
$$\lno{&u(x)(\omega)=\pi_{\ell(x,\omega)}(x)(\omega)&(4.10)\cr
&v(x)(\omega)=x(\omega)-u(x)(\omega).&(4.11)\cr}$$

\n{\bf Step 3.}\quad We fix $x$ and we show that $\Vert v(x)\Vert_1\le
K\theta_i(T)$.  
We define 
$$\lno{m(x,\omega)=\inf\left\{j\ge i-1;\vert x(\omega)-\pi_{j+1}(x)(\omega)
\vert > {r^{-j-1}\over2}\right\}&&(4.12)\cr}$$
(when the set on the right is empty, we set $m(x,\omega)=\infty$).

Since we assume $T\subset r^{-i}B_{\infty}/4$, we have 
$$\vert x(\omega)-\pi_i(x)(\omega)\vert\le {r^{-i}\over 2}$$
and this shows that $m(x,\omega)\ge i$.  The definition of $m(x,\omega)$ 
thereby implies that 
$$\vert x(\omega)-\pi_{m(x,\omega)}(\omega)\vert\le {r^{-m(x,\omega)}\over 
2}.$$

For $j<m(x,\omega)$, we have 
$$\eqalign{\vert\pi_j(x)(\omega)-\pi_{j+1}(x)(\omega)\vert&\le \vert
x(\omega)-\pi_j(x)(\omega)\vert+\vert x(\omega)-\pi_{j+1}(x)(\omega)
\vert\cr
&\le {r^{-j}\over 2}+{r^{-j-1}\over 2}<r^{-j},\cr}$$
so that $j<\ell(x,\omega)$.  This shows that $m(x,\omega)\le \ell(x,
\omega)$.

We show now that 
$$\lno{\vert v(x)(\omega)\vert\le 2r^{-m(x,\omega)}.&&(4.13)\cr}$$

Indeed
$$\eqalign{\vert v(x)(\omega)\vert&=\vert
x(\omega)-\pi_{\ell(x,\omega)}(x)(\omega)\vert\cr
&\le \vert x(\omega)-\pi_{m(x,\omega)}(x)(\omega)\vert +\sum_{m(x,\omega)
\le \ell<\ell(x,\omega)}\vert\pi_{\ell}(x)(\omega)-\pi_{\ell+1}(x)(\omega)
\vert\cr
&\le {1\over 2}r^{-m(x,\omega)}+\sum_{\ell\ge m(x,\omega)}r^{-\ell}\le
2r^{-m(x,\omega)}.\cr}$$

For $j\ge i-1$, we set $A_j=\{\omega;m(x,\omega)=j\}$.  By definition of
$m(x,\omega)$, we have 
$$\omega\in A_j\Rightarrow \vert
x(\omega)-\pi_{j+1}(x)(\omega)\vert>{r^{-j-1}\over 2}.$$

Since $x,\pi_{j+1}(x)$ both belong to $C_{j+1}(x)$, we use (4.7) with
$a=r^{-j-1}/2$ to see that 
$$\lambda(A_j)\le 4D_{j+1}(C_{j+1}(x)).$$
Thus 
$$\eqalign{&\int r^{-m(x,\omega)}d\lambda(\omega)=\sum_{j\ge
i-1}r^{-j}\lambda(A_j)\cr
&\le 4r\sum_{j\ge i-1}r^{-j-1}D_{j+1}(C_{j+1}(x))\le 8r\theta_i(T).\cr}$$
This proves that $\Vert v(x)\Vert_1\le Kr\theta_i(T)$ and finishes this 
step.

\n{\bf Step 4.}\quad We set $U=\{u(x);x\in T\}$, and we proceed to show 
that 
$$\lno{\gamma_1(U,\Vert\cdot\Vert_{\infty})\le K(\theta_i(T)+r^{-i}).
&&(4.14)\cr}$$

The main observation is that 
$$\lno{x,y\in C\in\cl C_j\Rightarrow \Vert u(x)-u(y)\Vert_{\infty}\le 4 
r^{-j}.&&(4.15)\cr}$$

Indeed, we have $\pi_{\ell}(x)=\pi_{\ell}(y)$ for $\ell\le j$.  Thus, the
definition of $\ell(x,\omega)$ shows that if either $\ell(x,\omega)<j$ or
$\ell(y,\omega)<j$ we have $\ell(x,\omega)=\ell(y,\omega)$, so that
$u(x)(\omega)=u(y)(\omega)$.  If $\ell(x,\omega)\ge j$, $\ell(y,\omega)\ge 
j$, we write 
$$\eqalign{\vert u(x)(\omega)-u(y)(\omega)\vert &\le \sum_{j\le \ell<
\ell(x,\omega)}\vert\pi_{\ell}(x)(\omega)-\pi_{\ell+1}(x)(\omega)\vert\cr
&\quad +\sum_{j\le \ell<
\ell(y,\omega)}\vert\pi_{\ell}(y)(\omega)-\pi_{\ell+1}(y)(\omega)\vert\cr
&\le 2\sum_{\ell\ge j}r^{-j}\le 4  r^{-j}.\cr}$$
This proves (4.15).  The result now follows simply by considering e.g. the 
probability that gives mass $2^{-j+i-1}\mu(C)$ to $u(y(C))$, for all $C\in
\cl C_j$, $j\ge i$, and by a routine computation.

\n{\bf Step 5.}\quad We show that 
$$\gamma_{1/2}(U)\le K(\theta_i(T)+r^{-i}+r^{-i+1}D_{i-1}(T)^{1/2}).$$

First, we observe that we can actually assume the following 
$$\text{If}~~C\supset D,~~~C\in\cl C_j,~~D\in \cl C_{j+1},~~\text{and~if}~~
y(C)\in D,~~\text{then}~~y(D)=y(C).$$

This implies in particular that
$$\lno{\pi_{\ell+1}(\pi_{\ell}(x))=\pi_{\ell}(x).&&(4.16)\cr}$$
for each $x,\ell$.

We fix $x$ and we estimate $\Vert u(x)-u(\pi_{\ell}(x))\Vert_2$.  We set
$G_{\ell}=\{\omega;\ell(x,\omega)=\ell\}$.  Since 
$$\omega\in G_{\ell}\Rightarrow
\vert\pi_{\ell}(x)(\omega)-\pi_{\ell+1}(x)(\omega)\vert > r^{-\ell},$$
it follows from (4.12), taking $a=r^{-\ell}$ that (since $\pi_{\ell}(x),
\pi_{\ell+1}(x)\in C_{\ell}(x)$) 
$$\lno{\lambda(G_{\ell})\le D_{\ell}(C_{\ell}(x)).&&(4.17)\cr}$$

We have seen in the proof of (4.15) that $u(x)(\omega)-u(\pi_j(x))(\omega)
=0$ unless $\ell(x,\omega)\ge j$.  We have, when $\ell(x,\omega)\ge j$ 
$$\lno{\vert u(x)(\omega)-u(\pi_j(x))(\omega)\vert&\le \sum_{j\le
\ell<\ell(x,\omega)}\vert\pi_{\ell}(x)(\omega)-\pi_{\ell+1}(x)(\omega)
\vert&(4.18)\cr
&\le \sum_{\ell\ge j}z_{\ell}(\omega)\cr}$$
where we define 
$$z_{\ell}(\omega)=\vert\pi_{\ell}(x)(\omega)-\pi_{\ell+1}(x)(\omega)\vert$$
if the right hand side is $\le r^{-\ell}$, and $z_{\ell}(\omega)=0$ 
otherwise.  In particular 
$$\Vert z_{\ell}\Vert^2_2\le \int \min (r^{-\ell},\vert
\pi_{\ell}(x)(\omega)-\pi_{\ell+1}(x)(\omega)\vert)^2d\lambda(\omega)\le
r^{-2\ell}D_{\ell}(C_{\ell}(x))$$
since $C_{\ell}(x)$ contains both $\pi_{\ell}(x)$ and $\pi_{\ell+1}(x)$.

Thus, by (4.17) and the triangle inequality we have 
$$\lno{\Vert u(x)-u(\pi_j(x))\Vert_2\le \sum_{\ell\ge
j}r^{-\ell}D_{\ell}(C_{\ell}(x))^{1/2}.&&(4.19)\cr}$$

\n Let us observe that this inequality holds in particular for $j=i-1$ so 
that 
$$\Vert u(x)-u(\pi_{i-1}(x))\Vert_2=\Vert u(x)-u(y(T))\Vert_2\le
r^{-i+1}D_{i-1}(T)^{1/2}+Kr^{-i/2}\theta_i(T)^{1/2}$$
since, for $\ell\ge i$, we have 
$$D_{\ell}(C_{\ell}(x))^{1/2}\le
r^{\ell/2}\theta_i(T)^{1/2}.$$

Thus using the inequality $\sqrt{ab}\le a+b$ we have shown the following:

\n (4.20)\quad The diameter of $U$ for $L^2$ is at most
$2 r^{-i+1}D_{i-1}(T)^{1/2}+K(\theta_i(T)+r^{-i})$.  

\n Consider now the measure $\nu$ that gives mass $2^{-j+i-1}\mu(C)$ to 
each point $u(y(C))$, $C\in \cl C_j$, $j\ge i$, and set 
$$w(x,j)=\log{1\over 2^{-j+i-1}\mu(C_j(x))}.$$
Since $u(\pi_j(x))$ has mass $\ge 2^{-j+i-1}\mu(C_j(x))$, setting $r(x)=
\Vert u(x)-u(\pi_i(x))\Vert_2$, we have, using the inequality $\sqrt{ab}
\le a+b$, as well as (4.19), 
$$\eqalign{I(x)=\int^{r(x)}_0\left(\log{1\over B_2(u(x),t)}
\right)^{1/2}dt&\le
\sum_{j\ge i}\Vert u(x)-u(\pi_j(x))\Vert_2\sqrt{w(x,j+1)}\cr
&\le 4\sum_{j\ge i}\left(\sum_{\ell\ge
j}r^{-\ell}(D_{\ell}(C_{\ell}(x)))^{1/2}\right)\sqrt{w(x,j+1)}\cr
&\le 4\sum_{j\ge i}\sum_{\ell\ge j}r^{-(\ell-j)/2}\sqrt{r^{-\ell}D_{\ell}
(C_{\ell}(x)) r^{-j}w(x,j+1)}\cr
&\le 8\sum_{j\ge i}\sum_{\ell\ge
j}r^{-(\ell-j)/2}(r^{-\ell}D_{\ell}(C_{\ell}(x))+r^{-j}w(x,j+1))\cr
&\le K(\sum_{\ell\ge i}r^{-\ell}D_{\ell}(C_{\ell}(x))+\sum_{j\ge i}r^{-j}
w(x,j+1))\cr
&\le K(r\theta_i(T)+r^{-i}).\cr}$$
The conclusion then follows from (4.20).\hfill\bx

The following is a kind of converse to Proposition 4.3 and proves that the
functional considered in this Proposition is indeed a sharp way to study 
the decompositions $T\subset U+uB_1$.  We now on assume that $\lambda$ is 
the counting measure on $\{1,\dots, N\}$.

\proclaim{Proposition 4.4}  Consider $i\in\Bbb Z$ and $U\subset L^2(\lambda
)$ with $\gamma_{1/2}(U)\le r^{-i}$.  Then 
$$\theta_i(U+r^{-i}B_1)\le K(r) r^{-i}.$$
\endproclaim

\demo{Proof}  {\bf Step 1.}\quad From the discussion prior to Proposition 
4.3 follows that $\theta_i(U)\le Kr^{-i}$.  Thus we can find an increasing 
sequence $(\cl A_j)_{j\ge i}$ of finite partitions of $U$, and a probability 
measure $\nu$ on $U$ such that 
$$\forall~x\in U,\quad \sum_{j\ge
i}r^{-j}\left(D_j(A_j(x))+\log{1\over\nu(A_j(x))}\right)\le Kr^{-i}.$$

\n{\bf Step 2.}\quad To each $y\in B_1$, we associate a sequence of 
integers $\ov p(y)=(p_{\ell}(y))_{\ell\ge 1}$ that satisfies the following 
properties.
$$\lno{&\lambda(\{r^{-\ell}\le \vert y\vert \le r^{-\ell+1}\})\le r^{{3
\over 2}p_{\ell}(y)}&(4.21)\cr
&\sum_{\ell\ge 1}r^{{3\over 2}p_{\ell}(y)-\ell}\le K(r)&(4.22)\cr
&\forall~\ell\ge 1,\quad \vert p_{\ell+1}(y)-p_{\ell}(y)\vert\le 1.
&(4.23)\cr}$$
To do this we denote by $q_{\ell}(y)$ the smallest integer for which the
left-hand side of (4.21) is $\le r^{{3\over 2}q_{\ell}(y)}$, and we set 
$$p_{\ell}(y)=\max\limits_{m\ge 1}(q_m(y)-\vert m-\ell\vert).$$
\enddemo

\n Thus (4.23) holds.  To prove (4.22), we simply observe that 
$$\eqalign{r^{{3\over 2}p_{\ell}(y)-\ell}&\le \sum_{m\ge 1}r^{{3\over
2}q_m(y)-{3\over 2}\vert m-\ell\vert-\ell}\cr
&\le \sum_{m\ge 1}r^{-\vert m-\ell\vert/2}r^{{3\over 2}q_m(y)-m}\cr}$$
and we invert the summation signs.

\n{\bf Step 3.}\quad For $j\ge i$, we consider the family $S_j$ of 
sequences $\ov q$ of integers, $\ov q=(q_1,\dots, q_{j-i+1})$ that satisfy 
$$\sum_{1\le m\le j-i+1}r^{{3\over 2}q_m-m}\le K(r)$$
where $K(r)$ is the constant of (4.22).  

For each point $t\in U+r^{-i}B_1$, we choose once for all a decomposition 
$t=x(t)+r^{-i}y(t)$, where $x(t)\in U$, $y(t)\in B_1$.  (The choice is made 
arbitrarily among all possible decompositions.)   To $A\in\cl A_j$, $\ov q
\in S_j$, we associate the set $C(A,\ov q)$ that consists of all the points 
for which $x(t)\in A$ and $y=y(t)$ satisfies 
$$\lno{p_1(y)=q_1,\dots, p_{j-i+1}(y)=q_{j-i+1}.&&(4.24)\cr}$$
These sets $C(A,\ov q)$ form a finite partition of $U+r^{-i}B_1$, that we 
denote by $\cl C_j$.  The sequence $(\cl C_j)_{j\ge i}$ is increasing.

\n{\bf Step 4.}\quad We show that 
$$\lno{D_j(C(A,\ov q))\le K[D_j(A)+r^{{3\over 2}q_{j-i+1}}].&&(4.25)\cr}$$
Since $\varphi^{1/2}_j$ is a distance, it suffices to show that 
$$\varphi_j(r^{-i}y,0)\le 2r^{{3\over 2}q_{j-i+1}}$$
whenever (4.24) holds.

To simplify notations, we write $p_m$ rather than $p_m(y)$.  We have 
$$\lno{\varphi_j(r^{-i}y,0)&=\int \min(r^{2j-2i}y^2,1)d\lambda&(4.26)\cr
&\le \sum_{\ell\ge 1}r^{-2\ell+2}\lambda(\{r^{2j-2i}y^2\ge r^{-2\ell}\}).
\cr}$$
Now, 
$$\eqalign{\lambda(\{r^{2j-2i}y^2\ge r^{-2\ell}\})&=\lambda(\{y\ge
r^{-\ell+i-j}\})\cr
&\le r^{{3\over 2}p_{\ell+j-i}}.\cr}$$
Thus, by (4.26)
$$\eqalign{\varphi_j(r^{-j}y,0)&\le r^2\sum_{\ell\ge 1}r^{-2\ell+{3\over
2}p_{\ell+j-i}}\cr
&=r^{2j-2i+2}\sum_{\ell\ge 1}r^{{3\over 2}p_{\ell+j-i}-2(\ell+j-i)}.\cr}$$
It follows from (4.23) that this sum is at most twice its first term, so 
that 
$$\varphi_j(r^{-j}y,0)\le 2r^{{3\over 2}p_{j-i+1}}=2r^{{3\over 2}
q_{j-i+1}}$$
since $p_{j-i+1}=p_{j-i+1}(y)=q_{j-i+1}$.

\n{\bf Step 5.}\quad The definition of $(\cl C_j)$ shows that if $x = x(t), 
y = y(t)$, then, for $j\ge i$, we have $C_j(t)=C(A_j(x),\ov q^j(y))$, where 
$$\ov q^j(y)=(p_1(y),\dots, p_{j-i+1}(y)).$$
Thus, by (4.25)
$$r^{-j}D_j(C_j(t))\le K[r^{-j}D_j(A_j(x))+r^{{3\over 2}p_{j-i+1}(y)-j}].$$
Combining with (4.22) yields 
$$\sum_{j\ge i}r^{-j}D_j(C_j(t))\le K(r)(\theta_i(U)+r^{-i}).$$

\n{\bf Step 6.}  We observe that $\card S_j\le (K(r)(j-i+1))^{j-i}$.  There 
is a probability $\mu$ on $A+r^{-i}B_1$ that gives mass at least 
$$\nu(C)2^{-j+i-1}(\card S_j)^{-1}$$
to each set $C(A,\ov q)$, $\ov q\in S_j$, $A\in \cl A_j$.  The fact that 
$$\sum_{j\ge i}r^{-j}\log{1\over\mu(C_j(t))}\le K(r)(\theta_i(U)+r^{-i})$$
then follows by a routine computation.\hfill\bx 

The measure of the size of a set of functions by the quantity $\theta_i(T)$
associated to functionals $\varphi_j$ given by (4.9) seems to be the
correct way to capture decompositions $T\subset U+uB_1$, in Propositions 
4.3 and 4.4.  However, in order to prove Theorem 1.3 it will be easier to 
use a cruder tool, that is exactly adapted to the study of weaker 
decompositions $T\subset U+uB_{p,\infty}$, where 
$$B_{p,\infty}=\{f\in L^0(\lambda);\sup\limits_{t\ge 0}t^p\lambda(\{\vert
f\vert\ge t\})\le 1\}.$$

We consider $1\le p<2$, and we define $\gamma$ by $\gamma(2-p)=1$.  We set 
$$d_i(f,g)=(\int \min((f-g)^2,r^{-4\gamma i})d\lambda)^{1/2}.$$

\proclaim{Proposition 4.5}  Consider a subset $T$ of $L^2(\lambda)$.  
Assume that $T\subset {1\over 4}B_{\infty}$.  Assume that there is an 
increasing sequence of finite partitions $(\cl C_j)_{j\ge 0}$ of $T$ and a 
probability measure $\mu$ on $T$ such that the following holds for a 
certain number $S\ge 1$.
$$\lno{&\text{Each set $C\in\cl C_j$ is of diameter $\le r^{-j}$ for 
$d_j$}. &(4.27)\cr
&\forall~x\in T,\quad \sum_{j\ge
0}r^{-j}\left(\log{1\over\mu(C_j(x))}\right)^{1/2}\le S.&(4.28)\cr}$$

Then we can find a set $U$ with $\gamma_{1/2}(U)\le K(p)S$ such that 
$$T\subset U+K(p, r)B_{p,\infty}.$$
\endproclaim

\demo{Proof}  {\bf Step 1.}\quad For $C\in\cl C_j$, $j\ge 0$, we pick a 
point $y(C)\in C$.  For $x\in T$, we set $\pi_j(x)=y(C_j(x))$.  We define, 
for $x\in T$, $\omega\in\Omega$ 
$$\ell(x,\omega)=\inf\{j\ge 0;\vert \pi_j(x)(\omega)-\pi_{j+1}(x)(\omega)
\vert\ge r^{-2\gamma j}\}$$
and $\ell(x,\omega)=\infty$ when the set on the left is empty.  We set
$\pi_{\infty}(x)=x$, and we define 
$$\eqalign{&u(x)(\omega)=\pi_{\ell(x,\omega)}(x)(\omega)\cr
&v(x)=x-u(x).\cr}$$
\enddemo

\n{\bf Step 2.}\quad We set $U=\{u(x);x\in T\}$, and we proceed to prove 
that $\gamma_{1/2}(U)\le K(p)S$.  The basic fact is that, if $x,y\in C
\subset \cl C_k$, we have 
$$\Vert u(x)-u(y)\Vert_2\le K(p)r^{-k}.$$
Indeed, as in the proof of Proposition 4.3, we have $u(x)(\omega)=u(y)(
\omega)$ unless $\ell(x,\omega)\ge k$, $\ell(y,\omega)\ge k$.  Thus since 
$\pi_k(x) = \pi_k(y)$, it suffices to show that 
$$\lno{\Vert(u(x)-\pi_k(x))1_{\{\ell(x,\omega)\ge k\}}\Vert_2\le 2r^{-k}.
&&(4.29)\cr}$$
Now 
$$\lno{\vert u(x)(\omega)-\pi_k(x)(\omega)\vert 1_{\{\ell(x,\omega)\ge k
\}}&\le \sum_{\ell\ge k}\vert\pi_{\ell}(x)(\omega)-\pi_{\ell+1}(x)(\omega)
\vert 1_{\{\ell(x,\omega)>\ell\}}.&(4.30)\cr}$$
\enddemo

Since $\pi_{\ell}(x)$, $\pi_{\ell+1}(x)$ belong to $C_{\ell}(x)$, and since 
$\vert\pi_{\ell}(x)(\omega)-\pi_{\ell+1}(x)(\omega)\vert \le r^{-2\gamma
\ell}$ when $\ell<\ell(x,\omega)$, by (4.27) and the definition of $d_i$ we 
have
$$\Vert \pi_{\ell}(x)(\omega)-\pi_{\ell+1}(x)(\omega)
1_{\{\ell(x,\omega)>\ell\}}\Vert_2\le r^{-\ell}$$
so that (4.29) follows by the triangle inequality.

That $\gamma_{1/2}(U) \le K(p)S$ follows by the usual computation, putting 
mass $2^{-j-1}\mu(C)$ at each point $u(y(C))$.

\n{\bf Step 3.}\quad We prove that for each $x$, we have $\Vert
v(x)\Vert_{p,\infty}\le K(p)S$.  Consider 
$$m(x,\omega)=\sup\left\{j\ge 0;\vert x(\omega)-\pi_j(x)(\omega)\vert\le
{r^{-2\gamma j}\over 2}\right\}.$$
(Observe that the set on the right is not empty since $T\subset {1\over
4}B_{\infty}(\lambda)$.)  As in the proof of Proposition 2.3, we see that
$m(x,\omega)\le \ell(x,\omega)$.  Thus 
$$\lno{\vert v(x)(\omega)\vert&\le \vert
u(x)(\omega)-\pi_{m(x,\omega)}(x)\vert+\sum_{m(x,\omega)\le
\ell<\ell(x,\omega)}\vert\pi_{\ell}(x)(\omega)-\pi_{\ell+1}(x)(\omega)
\vert&(4.31)\cr
&\le K(p)r^{-2\gamma m(x,\omega)}.\cr}$$
If we set $H_k=\{m(x,\omega)=k\}$, we see that, by definition of $m(x,
\omega)$, we have 
$$\omega \in H_k \Rightarrow \vert x (\omega) - \pi_{k+1} (x) (\omega) 
\vert \ge {r^{-2 \gamma (k+1)} \over 2}$$
so that, using (4.7), 
$${r^{-4\gamma (k+1)}\over 4}\lambda(H_k)\le \int \min((x-\pi_{k+1}(x))^2,
r^{-4\gamma
(k+1)})d\lambda\le r^{-2(k+1)}$$
by the argument of (4.12) and since both $x$ and $\pi_{k+1}(x)$ belong to 
$C_{k+1} (x)$.  Thus 
$$\lno{\lambda(H_k)\le 4\  r^{4\gamma k-2k}=4 
r^{2(2\gamma-1)(k+1)}.&&(4.32)\cr}$$
Since, by the choice of $\gamma$ $(=(2-p)^{-1})$, we have
$2\gamma-1=p{\gamma}$, the result follows from (4.31) (4.32), since 
$$\eqalignno{\lambda(\{m(x,\omega)\le k\})\le \sum_{\ell\le
k}\lambda(H_{\ell})\le K(p,r)r^{2(2\gamma-1)k}.&&\square\cr}$$

The following shows that the method of Proposition 4.4 is indeed the 
correct approach to study the decompositions $T\subset U+uB_{p,\infty}$.

\proclaim{Proposition 4.4}  Consider a subset $T$ of $L^0(\lambda)$, and 
assume that $T\subset U+B_{p,\infty}$, where $\gamma_{1/2}(U)\le 1$.  Then 
one can find an increasing family of partitions $(\cl C_i)_{i\ge 0}$ of $C$ 
and a probability measure $\mu$ on $T$ such that, for each $C\in \cl C_i$, 
the diameter of $C$ for $d_i$ is $\le K(p)r^{-i}$, and that 
$$\forall~x\in T,\quad\sum_{i\ge
0}r^{-i}\left(\log{1\over\mu(C_i(x))}\right)^{1/2}\le K.$$
\endproclaim

\demo{Proof}  Using (4.3), it suffices to show that if a set $C$ has a
diameter for $L^2(\lambda)$ that is $\le 2  r^{-i}$, then $C+B_{p,\infty}$
has a diameter for $d_i$ that is $\le K(p)r^{-i}$.  Since $d_i$ is a 
distance, it suffices to observe that, for $f\in B_{p,\infty}$
$$\eqalignno{d_i(f,0)^2&=\int^{r^{-2\gamma i}}_0\lambda(\{f\ge t\})d(t^2)
\le \int^{r^{-2\gamma i}}_0t^{-p}2t\,dt\cr
&={2\over 2-p}r^{-2\gamma i\cdot (2-p)}={2\over 2-p}r^{-2i}.&\square\cr}$$
\enddemo
\bigskip

\subhead 5.~~Construction of majorizing measures in the two-parameters
situation\endsubhead

We consider a set $T$, and a family $\varphi_j$ of functions from $T\times 
T$ to $\Bbb R^+$ $(j\in \Bbb Z)$.  We assume $\varphi_j(s,t)=\varphi_j(t,
s)$, and 
$$\lno{\forall~s,t,u\in T,\quad  \varphi_j(s,t)\le
2(\varphi_j(s,u)+\varphi_j(u,t)).&&(5.1)\cr}$$
For a subset $S$ of $T$, we set $D_j(S)=\sup\{\varphi_j(s,t);s,t\in S\}$.  
Given $t\in T$, $a>0$, we set 
$$B_j(t,a)=\{s\in T;\varphi_j(t,s)\le a\}.$$
Thus, by (5.1) we have 
$$\lno{D_j(B_j(t,a))\le 4a.&&(5.2)\cr}$$

We consider a number $r\ge 2$, and for simplicity we assume that $r$ is a 
power of $2$ $(r=2^{\tau},\tau\in\Bbb N)$.  We make the crucial assumption 
that, for some $\delta>0$, 
$$\lno{\forall~s,t\in T,~~\forall~j\in\Bbb Z,\quad \varphi_{j+1}(s,t)\ge
r^{1+\delta}\varphi_j(s,t).&&(5.3)\cr}$$
The functions $\varphi_i$ given by (4.10) {\it do not} satisfy this 
condition.

We assume that to each subset $S$ of $T$ is associated a number $F(S)\ge 0$. 
If $S\subset S'$, we assume $F(S)\le F(S')$.  We consider an increasing 
family $(\cl A_j)_{j\ge i}$ of finite partitions of $T$.  We assume the 
following condition, where $\alpha,\beta$ are $>0$ (this condition is a 
substitute for (2.4)).

\n (5.4)\quad Consider $j\ge i$, and consider $p\ge \tau-1$.  
Consider a subset $C$ of $T$, and assume that for a certain $D\in\cl A_j$, 
we have $C\subset D$.  Assume that $D_{j-1}( C)\le 2^{p-\tau+3}$.  Set 
$N=2^{2^p}$, and consider points $t_1,\dots, t_N$ of $C$, such that 
$$\ell,\ell'\le N,~~\ell\not= \ell'\Rightarrow \varphi_j(t_{\ell},t_{\ell'}
)\ge 2^p.$$
Consider for each $\ell\le N$, a subset $A_{\ell}$ of
$C\cap B_j(t_{\ell},\alpha 2^p)$.  Then 
$$F(\bigcup\limits_{\ell\le N}A_{\ell})\ge \beta r^{-j}2^p+
\min\limits_{\ell\le N}F(A_{\ell}).$$

\proclaim{Theorem 5.1}  Suppose, with the notations above that $\alpha
r^{\delta}\ge 4$, and consider a probability $\nu$ on $T$, and an 
increasing sequence $(\cl A_j)_{j\ge i}$ of finite partitions of $T$.  
Then there exists a probability $\mu$ on $T$ and an increasing
sequence of finite partitions $(\cl C_j)_{j\ge i}$ on $T$ such that 
$$\lno{&\forall~x\in T,\quad\sum_{j\ge
i}r^{-j}\left(D_j(C_j(x))+\log{1\over\mu(C_j(x))}\right)&(5.5)\cr
&\le K\left({1\over \beta}F(T)+r^{-i}(1+D_{i-1}(T))+\sup\limits_{y\in T}
\sum_{j\ge i}r^{-j}\log{1\over \nu(A_j(y))}\right).\cr}$$
\endproclaim

We first present the basic construction.  (This construction will then be
iterated to prove Theorem 5.1.)

\proclaim{Proposition 5.2}  Consider $j\ge i$, a subset $C$ of $T$.  Assume 
that $C$ is contained in a set belonging to $\cl A_j$.  Assume that 
$D_{j-1}(C)\le 2^{n+2}$, and that we are given a number $a(C)$ that 
satisfies the following two properties 
$$\lno{&F(C)-\beta r^{-j}2^n\le a(C)\le F(C)&(5.6)\cr
&\forall~t\in C,\quad F(C\cap B_j(t,\alpha 2^{n-1}))\le
a(C)+\varepsilon_j&(5.7)\cr}$$
where $\varepsilon_j=\beta r^{-i}2^{-j+i}$.

Set $n'=n+\tau-1$.  Then, for $s\ge n'$, $\ell\le N_s=2^{2^s}$, we can find 
sets $V(s,\ell)$ and numbers $a(V(s,\ell))$ that satisfy the following 
conditions
$$\lno{&\text{The sets $V(s,\ell)$, $s\ge n'$, $\ell\le N_s$ form a 
partition of $C$}&(5.8)\cr
&D_{j+1}(V(s,\ell))\le 2^{s+2}&(5.9)\cr
&F(V(s,\ell))-\beta r^{-j-1}2^s\le a(V(s,\ell))\le F(V(s,\ell))&(5.10)\cr
&\forall~t\in V(s,\ell),\quad F(V(s,\ell)\cap B_{j+1}(t,\alpha 2^{s-1}))
\le a(V(s,\ell))+\varepsilon_{j+1}&(5.11)\cr
&F(V(s,\ell))+a(V(s,\ell))+{\beta\over 4}r^{-j-1}2^s\le F(C)+a(C)+{\beta
\over 8}r^{-j}2^n+\varepsilon_j.&(5.12)\cr}$$
\endproclaim

\subhead Comments\endsubhead 1)\quad Conditions (5.10), (5.11) express that
$a(V(s,\ell))$ is to $V(s,\ell)$ what $a(C)$ is to $C$.

\n 2)\quad The reader observes the different coefficients of $\beta$ in 
(5.12) so that summation of such relations does yield information.

\demo{Proof}  {\bf Step 1.}\quad {\it Construction.}\quad Starting with 
$s=n'$, we construct points $t(s,\ell)$ of $C$, $\ell\le N_s$, that satisfy 
the following conditions.

\n (5.13)\quad Denote by $H(s,\ell)$ the union of the sets
$B_{j+1}(t(s',\ell'),2^{s'})$ for either $s'<s$ or $s=s'$, $\ell'<\ell$.  
Then 
$$t(s,\ell)\not\in H(s,\ell).$$
$$\lno{F(C\cap B_{j+1}(t(s,\ell),\alpha 2^s))\ge \sup\{F(C\cap B_{j+1}(t,
\alpha 2^s));t\in C\bs H(s,\ell)\}-\varepsilon_{j+2}.&&(5.14)\cr}$$
The construction is immediate.  It continues as long as possible.  We set 
$$\eqalign{&W(s,\ell)=C\cap B_{j+1}(t(s,\ell),2^s)\cr
&V(s,\ell)=W(s,\ell)\bs \bigcup W(s',\ell')\cr}$$
where the union is over all the choices of $s'<s$ or $s'=s$, $\ell'<\ell$.  
Itis obvious that the sets $V(s,\ell)$ form a partition of $C$.  Also, by 
(5.2) we have $D_{j+1}(V(s,\ell))\le 2^{s+2}$. 

\n{\bf Step 2.}\quad If $s=n'=n+\tau-1$, we set $a(V(s,\ell))=F(V(s,\ell))$.  

If $s\ge n'+1=n+\tau$, we set 
$$\lno{a(V(s,\ell))=\min(F(V(s,\ell)),F(C)-\beta r^{-j-1}2^{s-1})
&&(5.15)\cr}$$
(so that (5.10) holds)~and we prove (5.11).  It suffices to consider the 
case $s\ge n'+1$, and to prove that 
$$\forall t \in V (s, t), F(V(s, \ell) \cap B_{j+1}(t, \alpha 2^{s-1})) 
\le F(C) - \beta r^{-j-1} 2^{s-1} + \varepsilon_{j+1}.$$  
\enddemo

First, we observe that, by construction of $t(s-1,N_{s-1})$, (condition 
(5.14)) we have 
$$\lno{\forall~t\in C\bs &H(s-1,N_{s-1}),&(5.16)\cr
&F(C\cap B_{j+1}(t,\alpha 2^{s-1})) \le F(C\cap B_{j+1}
(t(s-1,N_{s-1}),\alpha 2^{s-1}))+\varepsilon_{j+2}.\cr}$$
By (5.13), we have, for $\ell,\ell'\le N_{s-1}$, 
$$\ell\not= \ell'\Rightarrow \varphi_{j+1}(t(s-1,\ell),t(s-1,\ell'))\ge 
2^{s-1}.$$

\n We set 
$$A_{\ell}=C\cap B_{j+1}(t(s-1,\ell),\alpha 2^{s-1}).$$
We see that we can use (5.4) (with $j+1$ rather than $j$, $s-1$ rather than 
$p$).  Indeed, $D_j(C)\le 2^{n+2}\le 2^{s-\tau+2}$ since $s\ge n+\tau$.

Thus, we have 
$$\lno{F(C)\ge \beta r^{-j-1}2^{s-1}+\min\limits_{\ell\le
N_{s-1}}F(A_{\ell}).&&(5.17)\cr}$$
On the other hand, using (5.14) again, we have 
$$\ell'<\ell\Rightarrow F(C\cap B_{j+1}(t(s-1,\ell),\alpha 2^{s-1}))
\le F(C\cap B_{j+1}(t(s-1,\ell'),\alpha 2^{s-1}))+\varepsilon_{j+2}.$$
Thus 
$$\lno{\varepsilon_{j+2}+\min\limits_{\ell\le N_{s-1}}F(A_{\ell})\ge F(C
\cap B_{j+1}(t(s-1,N_{s-1}),\alpha 2^{s-1}))&&(5.18)\cr}$$
and combining (5.16) to (5.18) we get the result.

\n{\bf Step 3.}\quad We show that if $s\le n+\tau+1$ we have 
$$\lno{F(V(s,\ell))\le a(C)+\varepsilon_{j}.&&(5.19)\cr}$$
Indeed we have 
$$V(s,\ell)\subset C\cap B_{j+1}(t(s,\ell),2^s).$$
Now, by (5.3), we have 
$$B_{j+1}(t(s,\ell),2^s)\subset B_j(t(s,\ell),r^{-1-\delta}2^s).$$
Since $\alpha r^{\delta}\ge 4$ and $r=2^\tau$ we have $r^{-1-\delta}2^s\le 
\alpha 2^{s-\tau-2}\le \alpha 2^{n-1}$, so the result
follows from (5.7).

\n{\bf Step 4.}\quad We prove (5.12).  For this, we must distinguish cases.

\n{\bf Case 1.}\quad $s=n'=n+\tau-1$.  We note that by (5.10) we have 
$a(V(s, \ell)) \le F(V(s, \ell)) \le F(C)$.  By (5.19), we have 
$$F(V(s,\ell))+a(V(s,\ell))\le F(C)+a(C)+\varepsilon_{j}.$$
Since $r^{-j-1}2^s={1\over 2}r^{-j}2^n$, (5.12) follows.

\n{\bf Case 2.}\quad $n+\tau\le s\le n+\tau+1$.  By definition of
$a(V(s,\ell))$, we have
$$\lno{a(V(s,\ell))\le F(C)-\beta r^{-j-1}2^{s-1}.&&(5.20)\cr}$$

Combining with (5.19) we get 
$$F(V(s,\ell))+a(V(s,\ell))\le F(C)+a(C)-\beta r^{-j-1}2^{s-1}+
\varepsilon_j$$
from which (5.12) follows.

\n{\bf Case 3.}\quad $s\ge n+\tau+2$.  By (5.6), we have 
$$a(C)\ge F(C)-\beta r^{-j}2^n.$$
Since $F(V(s,\ell))\le F(C)\le a(C)+\beta r^{-j}2^n$, we then have 
$$a(V(s,\ell))+F(V(s,\ell))\le a(C)+F(C)+\beta r^{-j}2^n-\beta
r^{-j-1}2^{s-1}.$$
Since $s\ge n+\tau+2$, we have $r^{-j-1}2^{s-1}\ge 2r^{-j}2^n$, so that 
$$\beta r^{-j}2^n-\beta r^{-j-1}2^{s-1}\le -{\beta\over 4}r^{-j-1}2^s$$
from which (5.12) follows.  This completes the proof of Proposition 
5.2.\hfill\bx

\demo{Proof of Theorem 5.1}  The construction of the family $\cl C_j$
goes by induction over $j$.  Together with each element $D$ of $\cl C_j$,
we will also construct an index $n(D)$ and a positive number $c(D)$, in
such a way that the following conditions hold:
$$\lno{&D_j(D)\le 2^{n(D)+2}&(5.21)\cr
&F(D)-\beta r^{-j}2^{n(D)}\le c(D)\le F(D)&(5.22)\cr
&\forall~t\in D,\quad F(D\cap B_j(t,\alpha 2^{n(t)-1}))\le
c(D)+\varepsilon_j.&(5.23)\cr}$$
\enddemo

We start the construction with $\cl C_{i-1}=\{T\}$, $c(T)=F(T)$ and for
$n(T)$ the smallest integer such that $D_{i-1}(T)\le 2^{n(T)+2}$.

Assume now that we have constructed $\cl C_j$.  We show how to partition
a given element $D$ of $\cl C_j$.  First, we break $D$ into the pieces
$D\cap A$, $A\in \cl A_j$.  We fix $A$, and we show how to partition
$C=D\cap A$.  We set 
$$a(C)=\min (c(D),F(C)).$$

Thus, setting $n=n(D)$, we have by (5.22)
$$\lno{&F(C)-\beta r^{-j}2^n\le a(C)\le F(C)&(5.24)\cr
&\forall~t\in C,\quad F(C\cap B_j(t,\alpha 2^{n-1}))\le
a(D)+\varepsilon_j.&(5.25)\cr}$$
We are then in a position to apply Proposition 5.2.  The partition
$(V(s,\ell))$, $s\ge n+\tau-1$, $\ell\le N_s$ that we obtained is the
partition we want.  We set 
$$n(V(s,\ell))=s;\quad c(V(s,\ell))=a(V(s,\ell)).$$
This completes the construction of $\cl C_{j+1}$.  We observe that, by
construction 
$$\lno{\forall~D\in\cl C_j,\quad \forall~A\in\cl A_j,\quad \card
\{C\in\cl C_{j+1}; C\subseteq D\cap A, n(C)=n\}\le N_n=2^{2^n}.
&&(5.26)\cr}$$
By (5.9) to (5.11), (5.21) and (5.22) will hold for any element $D'$ of
$\cl C_{j+1}$ (so that the construction can continue).  We rewrite (5.12)
as 
$$\eqalign{F(D')+c(D')+{\beta\over 4}r^{-j-1}2^{n(D')}&\le
F(C)+a(C)+{\beta\over 8}r^{-j}2^{n(D)}+\varepsilon_j\cr
&\le F(D)+c(D)+{\beta\over 8}r^{-j}2^{n(D)}+\varepsilon_j.\cr}$$

Thus, for any $x\in T$, we have 
$$\eqalign{F(C_{j+1}(x))&+c(C_{j+1}(x))+{\beta\over
4}r^{-j-1}2^{n(C_{j+1}(x))}\cr
&\le F(C_j(x))+c(C_j(x))+{\beta\over
8}r^{-j}2^{n(C_j(x))}+\varepsilon_j.\cr}$$
We sum these inequalities for $j\ge i-1$.  We get 
$$\lno{{\beta\over 8}\sum_{j\ge i}r^{-j}2^{n(C_j(x))}\le 2F(T)
+{\beta\over 8}r^{-i-1}2^{n(T)}+ 4\beta r^{-i}&&(5.27)\cr}$$
so that 
$$\lno{\sum_{j\ge i}r^{-j}2^{n(C_j(x))}\le {16\over
\beta}F(T)+Kr^{-i-1}(1+D_{i-1}(T)).&&(5.28)\cr}$$

We now construct the measure $\mu$.  First, by induction over $j$ we
construct weights $w(D)$ for $D\in \cl C_j$, such that 
$$\lno{\sum_{D\in\cl C_j}w(D)\le 2^{i-j-2}.&&(5.29)\cr}$$
To do this, if $D'\in \cl C_{j+1}$ is contained in $C\cap A$ $(D\in\cl
C_j,A\in\cl A_j)$, we set 
$$\lno{w(D')={1\over 4}\nu(A)w(D)2^{-2^{n(D')+1}}.&&(5.30)\cr}$$
It follows from (5.26) that 
$$\sum_{D'\subset D\cap A}2^{-2^{n(D')+1}}\le \sum_{n\ge
0}2^{2^n}2^{-2^{n+1}}=\sum_{n\ge 0}2^{-2^n}\le \sum_{n\ge 0}2^{-n}=2$$
so that (5.29) at rank $j+1$ follows from (5.30) and the fact that 
$\mu$ is a probability.

It follows from (5.29) that there exists a probability $\mu$ on $T$ that
gives mass $\ge w(D)$ to each $D\in \cl D_j$.  From (5.30) we have 
$$\log{1\over\mu(C_{j+1}(x))}\le \log{1\over w(C_{j+1}(x))}\le \log
4+2^{n(C_{j+1}(x))}+\log{1\over \nu(A(x))}+\log{1\over w(C_j(x))}.$$
Using (5.28), (5.22), the conclusion follows easily.\hfill\bx
\medskip

Before we prove Theorem 1.2, we should mention that Theorem 5.1 provides
a simpler proof of Theorem 5.1 of [T5].  To see this, we take $\cl
A_j=\{T\}$ for all $j$, and we define for $F(S)$ the ``size'' of the 
largest tree (in the precise sense of Theorem 5.1 of [T4]) that is 
contained in $S$.  That (5.1) holds is obvious from the definition of the 
``size'' of a tree.

An essential ingredient to the proof of Theorem 1.2 is the following, that 
is a weakening of Corollary 2.7 of [T5].

\proclaim{Proposition 5.3}  Consider $t_1,\dots, t_N\in\Bbb R^M$, and
$a,b>0$.  Assume 
$$\lno{&\forall~\ell,\ell'\le N,\quad \Vert
t_{\ell}-t_{\ell'}\Vert_{\infty}\le a&(5.31)\cr
&\forall~\ell\not= \ell',\quad\Vert t_{\ell}-t_{\ell'}\Vert_2\ge
b.&(5.32)\cr}$$
Consider $\sigma>0$, and sets $A_{\ell}\subset B_2(t_{\ell},\sigma)$ for
$\ell \le N$.  Then 
$$b(\bigcup\limits_{\ell\le N}A_{\ell})\ge {1\over K}\min\left(b\sqrt{\log
N},{b^2\over a}\right)+\min\limits_{\ell\le
N}b(A_{\ell})-K\sigma\sqrt{\log N}.$$
\endproclaim

\proclaim{Corollary 5.4}  If $\sigma\le b/K_1$, $a\le 2b/\sqrt{\log N}$,
we have 
$$\lno{b(\bigcup\limits_{\ell\le N}A_{\ell})\ge {1\over K}b\sqrt{\log
N}+\min\limits_{\ell\le N} b(A_{\ell}).&&(5.33)\cr}$$
\endproclaim

\demo{Proof of Theorem 1.2}

\n{\bf Step 1.}\quad We choose $\alpha=1/K^2_1$, where $K_1$ occurs in
Corollary 5.4, and we chose for $r$ the smallest power of $2$ such that
$r\alpha\ge 4$.  We define $i$ as the smallest for which $r^{-i}$ is
larger than the diameter of $T$ for $\Vert \cdot\Vert_{\infty}$.  First, 
we find an increasing sequence $(\cl A_j)_{j\ge i}$ of finite
partitions of $T$, and a probability measure $\nu$ on $T$ such that 
$$\forall~x\in T,\quad\sum_{j\ge i}r^{-j}\log{1\over\nu(A_j(x))}\le
K\gamma_1(T,\Vert\cdot\Vert_{\infty})$$
and that the diameter for $\Vert\cdot\Vert_{\infty}$ of each $A\in 
\cl A_j$ is at most $2r^{-j}$.  We set $\varphi_j(s,t)=2^{2j}\Vert
s-t\Vert^2_2$, so that (5.3) holds with $\delta=1$.
\enddemo

We now prove that (5.4) holds for a certain $\beta>0$, when $F(S)=b(S)$. 
Since it is assumed in (5.4) that $C$ is contained in a set of $\cl A_j$, 
we see that (5.31) holds for $a=2r^{-j}$.  The definition of
$\varphi_j$ shows that (5.32) holds for $b=r^{-j}2^{p/2}$.  Since
$A_{\ell}$ is contained in $B_j(T_{\ell},\alpha 2^p)$, we see that the
number $\sigma$ of Proposition 5.3 can be taken equal to
$\sqrt{\alpha}\,r^{-j}2^{p/2}$, so that, by the choice of $\alpha$, we
have $\sigma\le b/K_1$.  Since $\sqrt{\log N}=2^{p/2}\sqrt {\log 2}$, we 
have $a\le 2b/\sqrt{\log N}$, so that the result follows from (5.33). 

Thus, we can use Theorem 5.1.  The right-hand side of (5.5) is at most 
$$K(b(T)+r^{-i}+r^i\Delta^2+\gamma_1(T,\Vert \cdot\Vert_{\infty}))$$
where $\Delta$ is the $\ell_2$ diameter of $T$.  Now, $r^{-i}\le
K\gamma_1(T,\Vert\cdot\Vert_{\infty})$, and since the $\ell_2$ diameter is
less than the $\ell_{\infty}$ diameter, this is at most 
$$\eqalignno{K(b(T)+\gamma_1(T,\Vert\cdot\Vert_{\infty})).&&\square\cr}$$
\bigskip

We now turn to the proof of Theorem 1.3.  Theorem 1.3 is a consequence of the 
following.

\proclaim {Proposition 5.5}  Consider vectors $(x_i)_{i \le M}$ is a 
Banach space $X$ of dimension $n$.  Then there is a subset $I$ of 
$\{1, \cdots , M\}$ with $\card I \le K n \log n$ such that either
$$\lno {E\Vert \sum_{i \not\in I, i \le M} g_i x_i \Vert \le {1 \over 2} 
E \Vert \sum_{i \le M} g_i x_i \Vert &&(5.34)\cr}$$
or 
$$\lno {E \Vert \sum_{i \not\in I, i \le M}  g_i x_i \Vert \le 
K E \Vert \sum_{i \le M} \varepsilon_i x_i \Vert .&& (5.35)\cr}$$
\endproclaim

Indeed, to obtain Theorem 1.3, we simply iterate use of Proposition 5.5.  It  
is known that $E \Vert \sum_{i \le M} g_i x_i \Vert \le K \sqrt {\log 
(n+1)} E \Vert \sum_{i \le M} \varepsilon_i x_i \Vert.$

We start the proof of Proposition 5.5.  We set $X^\ast_1 = \{ x^\ast 
\in X^\ast; \Vert x^\ast \Vert \le 1\}$.  Consider a subset $I$ of 
$\{ 1, \cdots, M \}$ that will be chosen later, and set $J = \{ i \le M; 
i \not\in I\}$.  We set, for $x^\ast \in X^\ast$, 
$$\Vert x^\ast\Vert_\infty = \sup\limits_{i \in J} \vert x^\ast (x_i)\vert;~~ 
B_\infty = \{  x^\ast \in X^\ast; \Vert x^\ast \Vert_\infty \le 1\}.$$

We can then reformulate Theorem 1.2 as 
$$\lno {E\Vert \sum_{i \in J} g_i x_i \Vert \le K (E \Vert \sum_{i \in J} 
\varepsilon_i x_i \Vert + \gamma_1 (X^\ast_1, \Vert \cdot \Vert_\infty)).
&&(5.36)\cr}$$

For $x^\ast \in X^\ast$, we set
$$\lno {\Vert x^\ast \Vert^2_2 = \sum_{i \le M} x^\ast (x_i)^2&&(5.37)\cr}$$

(This is the $L^2$ norm associated to the gaussian random vector 
$\sum\limits_{i \le M} g_i x_i)$.  We set $B_2 = \{ x^\ast \in X^\ast; 
\Vert x^\ast 
\Vert_2 \le 1\}$.  The key to the proof is the following interpolation 
formulae

\proclaim {Lemma 5.6}  We have
$$\lno {\gamma_1 (X^\ast_1, \Vert \cdot \Vert_\infty) \le K 
\gamma_{1/2}(X^\ast_1, \Vert \cdot \Vert_2) \sup\limits_{\varepsilon > 0} 
\varepsilon \sqrt {\log N (B_2 , \varepsilon B_\infty)} &&(5.38)\cr}$$
\endproclaim

We will prove this later in order not to break the flow of the argument.  
We plug (5.38) into (5.36) remembering that $\gamma_{1/2}(X^\ast_1, \Vert 
\cdot \Vert_2) \le K E \Vert \sum_{i \le M} g_i x_i \Vert$ by (1.3).  
Thus we get 
$$\lno {E \Vert \sum_{i \in J} g_i x_i\Vert \le K E \Vert \sum_{i \in J} 
\varepsilon_i x_i \Vert + K \alpha E \Vert \sum_{i \le M} g_i x_i 
\Vert &&(5.39)\cr}$$
where 
$$\alpha = \sup\limits_{\varepsilon > 0} \varepsilon \sqrt {\log N(B_2, 
\varepsilon B_\infty)}$$
Thus, if we can arrange that $4K \alpha \le 1$, whenever
$E \Vert \sum\limits_{i \le M} g_i x_i \Vert \le 2E\Vert \sum\limits_{i 
\in J} g_i x_i \Vert$, (5.39) becomes 
$$E \Vert \sum_{i \in J} g_i x_i\Vert \le K E \Vert \sum\limits_{i \in J} 
\varepsilon_i x_i \Vert + {1 \over 2} E \Vert \sum\limits_{i \in J} 
g_i x_i \Vert$$
and this implies (5.35).

Before we study $\alpha$, we need some preliminaries.  The formulae (5.37) 
defines a semi-norm $\Vert \cdot \Vert_2$ on $X^\ast$.  By duality, this 
semi-norm defines a norm $\Vert \cdot \Vert_2$ on the linear span $H$ of 
the vectors $(x_i)_{i \le M}$.  The unit ball of that norm is the set of 
vectors $\sum\limits_{i \le M} \alpha _i x_i$ with $\sum\limits_{i \le M} 
\alpha^2_i \le 1$.  If we denote by $\nu$ the law of $\sum\limits_{i \le M} 
g_i x_i, 
\nu$ is a gaussian measure, $H$ is its reproducing kernel and $\Vert \cdot 
\Vert_2$ is the associated norm.  One way to reformulate (5.37) is to 
say that $\nu$ is the canonical gaussian measure on $H$, i.e. 
$$\l x^\ast , y^\ast\r = \int_H x^\ast(x) y^\ast(x) d \nu (x)$$

For a subset $A$ of $H$, we can measure its size $\ell (A)$ with respect 
to the canonical gaussian measure $\nu$ by
$$\lno {\ell (A) = \int_H \sup\limits_{x \in A} \l x, y\r d \nu (y)
&&(5.40)\cr}$$

Since $\nu$ is the law of $\sum\limits_{i \le M} x_i g_i$ we have
$$\lno {\ell(A) = E \sup\limits_{x \in A} \l x, \sum_{i \le M} x_i g_i \r
&&(5.41)\cr}$$
although this formulae will not be used in the present proof.  An important 
fact for the rest of this argument is that, denoting by $\conv A$ the 
balanced convex hull of a set $A$, for a sequence $(z_k)$ in $H$ we have 
$$\lno {\ell (\conv (z_k)) \le K \sup\limits_k \Vert z_k \Vert_2 
\sqrt {\log (k+1)}&&(5.42)\cr}$$
This results from a trivial computation; see e.g. the introduction of [T2]).

It follows from (5.37) that
$$\sum\limits_{i \le M} \Vert x_i \Vert^2_2 \le n $$
so that we can assume without loss of generality that $\Vert x_i \Vert \le 
\sqrt {n/i}$.  Consider a number $L$ to be adjusted later, and set 
$J = \{ i \le M ; i \ge L n \log n\}, C = \conv \{ x_i; i \in J\}$.  Then, by 
(5.42), we get 
$$\lno {\ell (C) \le K \sup\limits_{K \ge 1} \sqrt {{n \log (k+1) 
\over k + n L \log n}} \le K \sqrt { {\log L \over L}}&&(5.43)\cr}$$
as is easily seen by distinguishing the cases $k \le n L \log n$ and $k \ge 
n L \log n$.  What we need to remember from (5.43) is that $\ell (C)$ can 
be made arbitrarily  small taking $L$ large enough. 

To bound $\alpha$, we now simply apply the reverse Sudakov minoration 
as in [L-T], Chapter 4, (3.15), to obtain 
$$\alpha \le K \ell (B^\circ_\infty)$$
where $B^\circ_\infty$ is the polar of $B_\infty$, that, by the bipolar 
theorem, is exactly $C$.  The proof is complete.\hfill $\square$

\demo {Proof of Lemma 5.6}  Consider an increasing sequence  
$(\Cal A_n)_{n \ge n_0}$ of partitions of $T= X^\ast_1$, where each 
element $A$ of $\Cal A_n$ is of diameter (for $\Vert \cdot \Vert_2$) at
most $2^{-n}$, and a probability measure $\mu$ on $A$ such that
$$\lno {\sup\limits_{t \in T} \sum_{n \ge n_0} 2^{-n} \sqrt {\log {1 \over 
\mu (A_n(t))}} \le K \gamma_{1/2} (T, \Vert \cdot \Vert_2).&&(5.44)\cr}$$
There $A_n(t)$ denotes as usual the unique element of $\Cal A_n$ that 
contains $t$. Set 
$$\alpha = \sup\limits_{ \varepsilon > 0} \varepsilon 
 \sqrt {\log N(B_2, \varepsilon B_\infty)}.$$  
Consider the smallest 
integer $p_0$ with $2^{p_0} \alpha \ge 1$.  Given any $n \ge n_0$, any  
$A \in \Cal A_n$ and any  $p \ge p_0$ we can find a subset $F(A, p)$ 
such that 
$$\card (F(A, p)) \le \exp (2^{2p} \alpha^2)$$
such that each point of $A$ is within distance $2^{-n-p}$ of a point of 
$F(A, p)$.  For each choice of $A, p$, we put a mass 
$$\mu(A) 2^{-p+p_0 - n + n_0 - 2} e^{- \alpha^2 2^{2p}}$$
at each point of $F(A, p)$, for a total mass $\le 1$.

We now prove that the resulting measure witnesses (5.38).  Let us fix $t$ 
in $T$, 
and, for $n \ge n_0$, consider 
$$\varepsilon_n = {\alpha 2^{-n} \over \sqrt {\log {3 \over \mu 
(A_n(t))}}}$$
When $2^{-\ell} \le \varepsilon_n$, the ball $B_\infty (t, 2^{-\ell})$ of 
center $t$, of radius $2^{-\ell}$ 
for $\Vert \cdot \Vert_\infty$ satisfies by construction
$$\nu (B_\infty (t, 2^{-\ell})) \ge \mu (A_n (t))2^{-p - n + p_0 + n_0 - 2} 
e^{-\alpha^2 2^{2p}}$$
where $p = \ell - n$ (observe that $2^p \alpha \ge 1)$.  Thus 
$$\log {1 \over \nu (B_\infty (t, 2^{-\ell}))} \le \log {1 \over 
\mu (A_n (t))} + \ell - p_0 - n_0 + 2 + \alpha^2 2^{2(\ell - n)}$$
and thus 
$$S_n  =: \sum_{\varepsilon_{n-1} < 2^{-\ell} \le \varepsilon_n} 2^{-\ell} 
\log {1 \over \nu (B_\infty(t, 2^{-\ell}))} \le K \varepsilon_n 
\log {1 \over \mu (A_n(t))} + {\alpha^2 \over \varepsilon_{n+1}} 
2^{-2n} + a_n$$
where $a_n = \sum\limits_{\varepsilon_{n+1} \le 2^{-\ell} \le \varepsilon_n} 
2^{-\ell} ( \ell  - p_0 - n_0 + 2)$, so that 
$$S_n \le K \alpha [ 2^{-n} \sqrt {\log {2 \over \mu (A_n(t))}} 
+ 2^{-n-1} \sqrt {\log {2 \over \mu (A_{n+1}(t))}}]+ a_n.$$
Thus 
$$\sum_{n \ge n_0} S_n \le K \alpha \gamma_{1/2} (T, \Vert \cdot \Vert_2) + 
\sum_{n \ge n_0} a_n.$$

We now observe that $\varepsilon_{n_0}= \alpha 2^{-n_0}$, and 
that the diameter of $T$ for $\Vert \cdot \Vert_\infty$ is at most 
$ \varepsilon_{n_0}$ (since $N(B_2, \alpha B_\infty)=1)$.  Thus 
$$\int^\infty_0 \log {1 \over \mu (B_\infty(t, \varepsilon))} d \varepsilon 
\le K \sum_{n \ge n_0} S_n$$
Also, $\sum\limits_{n \ge n_0} a_n \le K \alpha 2^{-n_0} \le K \alpha 
\gamma_{1/2} (T, \Vert \cdot \Vert_2)$.  This completes the proof.
\hfill $\square$
\enddemo

The reader might have noticed that the argument of Proposition 5.5 shows 
that 
$$E \Vert \sum_{i \le M} \varepsilon_i g_i \Vert \le K E \Vert 
\sum_{i \le M} \varepsilon_i x_i \Vert$$
whenever $\ell (C) \le 1/K$, where $C = \{ x_i ; i \ge 1\}$.  (The quantity 
$\ell(C)$ is defined in the course of the proof of Proposition 5.5).  
However more is true. 

\proclaim {Theorem 5.7}  For vectors $(x_i)_{i \le M}$ in a Banach space, 
we have 
$$\lno {E \Vert \sum_{i \le M} g_i x_i \Vert \le K E \Vert \sum_{i \le M} 
\varepsilon_i x_i \Vert (1 + \ell (C))^3.&&(5.45)\cr}$$
\endproclaim

A positive solution to the Bernoulli problem would imply that (5.45) holds  
with a factor $(1+\ell (C))$ rather than $(1+\ell (C))^3$.  It seems, 
however, that the difficulties one faces in proving this (even after one 
has obtained the apparently optimal Lemma 5.8 below) are of the same nature 
as some of the difficulties one faces when studying the Bernoulli problem.  
On the other hand, we know how to do better than (5.45), and in particular 
how to prove 
$$\lno {E \Vert \sum_{i \le M} g_i x_i\Vert \le K E \Vert \sum_{i \le M} 
\varepsilon_i x_i \Vert (1 + \ell (C)) \log (2+ \ell (C)) &&(5.46)\cr}$$
The techniques to obtain this improved estimates are however not related 
to the other material of the present paper, but rather are variations on 
the ``tree extraction'' techniques of [T5].  Since, moreover, there is not 
much conceptual gain in proving the imperfect inequality (5.46) rather 
than the (slightly more imperfect) inequality (5.45), we will prove 
(5.45) only.

\demo {Proof of Theorem 5.7}  We set $\alpha = \ell (C)$.  A key estimate 
is as follows
\enddemo

\proclaim {Lemma 5.8}  Consider $y^\ast, y^\ast_1, \cdots , y^\ast_N$ in 
$X^\ast$.  Assume
$$\lno {\forall i, j \le N, \Vert y^\ast_i - y^\ast_j \Vert_2 \ge 1
&&(5.47)\cr}$$
$$\lno {E \sup\limits_{j \le N} y^\ast_j (\sum_{i \le M} \varepsilon_i 
x_i) \ge {\sqrt {\log N} \over K(1+\alpha)}&&(5.48)\cr}$$
\endproclaim

\demo {Proof}  Consider the map $W$ from $X^\ast$ to $\Bbb R^M$ that sends 
$y^\ast$ to $(y^\ast(x_i))_{i \le M}$.  Set $T = W(\{ y^\ast_j; j \le N\})$, 
so that the left-hand side of (5.48) is simply $b(T)$.  Denoting as usual 
by $B_2$ and $B_1$ the $\ell_2$ and the $\ell_1$ unit balls of $\Bbb R^M$, 
for numbers $\theta, L > 0$, consider $D = \theta B_2 + L B_1$.  Thus if 
$t \in D$, we can write
$$\lno {t_i = u_i + v_i; \sum_{ i \le M} u^2_i \le \theta^2, \sum_{i \le M} 
\vert v_i \vert \le L&&(5.49)\cr}$$
However, even if $t \in W(X^\ast)$ there is no reason why $(u_i)_{i \le M}$ 
or $(v_i)_{i \le M}$ should be of the same type.  This is why we moved to 
$\Bbb R^M$ rather than working in $X^\ast$.  

The key tool is the version of Sudakov minoration for Bernoulli processes 
proved in [T5] that asserts that
$$L \ge K b(T) \Rightarrow \theta \sqrt {\log N(T, D)} \le K b(T)$$ so that
$$\lno {b(T)\ge {1 \over K} \min (L , \theta \sqrt{\log N(T,D)})
&&(5.50)\cr}$$
Our task is now to find a lower bound for $N(T,D)$.  Consider $t \in T$, 
and $S = T\cap (t+D)$.  Set $R = \card S$. By (5.47) and (the usual) 
Sudakov minoration, we have
$$\lno {{1 \over K} \sqrt{\log R} \le G(S).&&(5.51)\cr}$$
To bound $\log R$, we now find an upper bound for the right-hand side 
of (5.51).  We follow the notation established during the proof of 
Proposition 5.5.  A basic observation is that for any $y^\ast$ in 
$X^\ast$ and any $x$ in $H$ we have
$$\lno {y^\ast (x)= \sum_{i \le M} y^\ast (x_i) \l x, x_i \r .&&(5.52)\cr}$$
This is a consequence of (5.37) and of the fact that the map $U$ from $H$ 
to $X^\ast$ given by $y^\ast(x) = \l y^\ast, U (x) \r$ defines an isometry 
from $(H, \Vert \cdot \Vert_2)$ into $(X^\ast, \Vert \cdot \Vert_2)$.

Consider now $s \in S$, so that $s - t \in D$.  Also, $s - t\in W(X^\ast)$.  
Thus there is $y^\ast$ in
 $X^\ast$ such that $i \le M$ we have $s_i- t_i = y^\ast(x_i)$, and thus, 
by (5.52)
$$\lno {\sum_{i \le M} (s_i - t_i) g_i &= y^\ast(\sum_{k \le M} g_k x_k) 
= \sum_{i \le M} y^\ast(x_i) \l x_i, \sum_{k \le M} g_k x_k \r &(5.53)\cr
& = \sum_{i \le M} (s_i - t_i) \l x_i, Z\r \cr}$$
where for simplicity we set $Z = \sum\limits_{k \le M} g_k x_k$.  We now 
appeal to (5.49), since $s-t \in D$, to write
$$s_i - t_i = u_i + v_i; \sum_{i \le M} u^2_i \le \theta^2 , \sum_{i \le M} 
\vert v_i \vert \le L.$$
and thus, by (5.53) we get
$$\sum_{i \le M} (s_i - t_i) g_i = \sum_{i \le M} u_i \l x_i, Z \r 
+ \sum_{i \le M} v_i \l x_i, Z \r$$
and thus 
$$\lno {\vert \sum_{i \le M} (s_i - t_i) g_i\vert &\le \vert \l 
\sum_{i \le M} u_i x_i, Z \r \vert + L \sup\limits_{i \le M}\vert \l x_i,
 Z \r&(5.55)\cr}$$ 
The r.v. 
$$Y_s = \l \sum_{i \le M} u_i x_i, Z\r= U(\sum_{i \le M} u_i x_i)(Z)$$ 
is gaussian, and by definition on the norm $\Vert \cdot \Vert_2$ on 
$X^\ast$,
$$E Y^2_s= \Vert U (\sum_{i \le M} u_i x_i) \Vert^2_2 = \Vert \sum_{i \le M} 
u_i x_i\Vert^2 \le \sum_{i \le M} u^2_i \le \theta^2$$
Now, we have, by (5.55)
$$\lno {G(S) = G(S-t) &= E \sup\limits_{s \in S} \sum_{i \le M} (s_i - t_i)
g_i \cr 
&\le E \sup\limits_{s \in S} \vert \sum_{i \le M} (s_i - t_i) g_i\vert\cr
&\le E \sup\limits_{s \in S} \vert Y_s \vert + L E \sup\limits_{i \le M} 
\vert\l x_i, Z \r \vert \cr}$$
Using (5.41) we see that
$$E \sup\limits_{i \le M} \vert \l x_i, Z \r \vert = \ell (C) = \alpha.$$
Thus, by a standard estimate
$$G(S) \le K \theta \sqrt {\log R} + L \alpha$$
and plugging back into (5.51) gives 
$$\lno {{1 \over K} \sqrt {\log R} \le K \theta \sqrt {\log R} + 
L\alpha.&&(5.56)\cr}$$
We now fix $\theta = 1/2 K^2, L = \sqrt {\log N}/2 K\alpha$, so that 
(5.56) implies
$$\sqrt {\log R} \le \sqrt {{ \log N \over 2}},$$
and hence $R \le \sqrt{N}$.  This shows that for this choices of $\theta, 
L, t + D$ contains at most 
$\sqrt {N}$ points of $T$; thus $N(T,D) \ge \sqrt {N}$, and plugging in 
(5.50) this proves the result.\hfill $\square$
\enddemo

Mimicking the proof of Corollary 2.7 of [T5] (and relying upon (5.7)) we 
obtain the following, where, for a subset $A$ of $X^\ast$, we set
$$b(A) = E \sup\limits_{x^\ast \in A} x^\ast (\sum_{i \le M} 
\varepsilon_i x_i).$$

\proclaim {Proposition 5.8}  Consider $y^\ast_1, \cdots , y^\ast_N$ in 
$X^\ast$.  Assume 
$$\forall \ell \not= \ell ', \Vert y^\ast_\ell - y^\ast_\ell 
\Vert_2 \ge b.$$
Consider $\sigma > 0$ and for $\ell \le N$ consider $A_\ell \subset
 B_2 (y^\ast, \sigma) = \{ z^\ast \in X^\ast; \Vert y^\ast - z^\ast 
\Vert_2 \le \sigma \}$.  Then
$$b(\bigcup\limits_{\ell \le N}A_\ell) \ge {b \over K(1+\alpha)} 
\sqrt {\log N} + \min\limits_{\ell \le N} b(A_\ell) - K \sigma 
\sqrt{ \log N}$$
\endproclaim

\proclaim {Corollary 5.9}  If $\sigma \le b/K (1+\alpha)$ we have 
$$b(\bigcup\limits_{\ell \le N}A_\ell) \ge \min\limits_{\ell \le N} 
b (A_\ell) + {b \over K(1+\alpha)} \sqrt {\log N}$$
\endproclaim

We now appeal to Theorem 2.1, with $p=1$.  We set $T = X^\ast_1$, provided 
with the distance induced by the norm $\Vert \cdot \Vert_2$, and we set
$$\varphi_k (x^\ast) = b(X^\ast_1)- b (X^\ast_1 \cap B_2 (x, r^{-k})).$$
It follows from Corollary 5.9 that (2.4) holds, provided $r = 
K(1+\alpha)$ and $$\theta (n) = {1 \over K r(1+ \alpha)} \sqrt {\log n}.$$

Now, we appeal to Theorem 2.2 and the remark that follows its proof.  We 
observe that an extra factor $r$ occurs when comparing the left-hand 
side of (2.19) with an integral such as the right-hand side of (1.1).  
This finishes the proof of Theorem 5.6. \hfill $\square$

\subhead 6.~~~Radmacher Cotype\endsubhead

To proof of Theorem 1.4 relies on a different version of the construction
of Section 1.  We consider a set $T$, such that on $T$ we have a sequence
$(d_j)_{j\ge 0}$ of distances.  We assume that this sequence is decreasing, 
i.e. $d_{j+1}(s,t)\le
d_j(s,t)$ for $s,t\in T$.  We denote by $B_j(x,a)$ the ball for
$d_j$.  We assume that for each $j\ge 0$, each subset
$S$ of $T$, we are given a quantity $F_j(S)$ that is increasing in $S$.  We
assume that the sequence $F_j$ of functionals is decreasing, i.e.
$F_{j+1}(S)\le F_j(S)$ for $S\subset T$, $j\ge 0$.  We assume that for
certain $\gamma>0$, $r\ge 4$, $K_2>0$, the following condition (that is a
substitute for (2.4)) holds.

\n (6.1)\quad Consider $j\ge 0$, $k\ge 0$, $t\in T$, $N\ge 2$, and points
$(t_{\ell})_{\ell\le N}$ in $B_j(t,r^{-k})$.

Assume that 
$$\ell\not= \ell'~~\text{implies}~~d_j(t_{\ell},t_{\ell'})\ge r^{-k-1}.$$
Consider sets $T_{\ell}\subset B_j(t_{\ell},r^{-k-2})$.  Then, whenever 
$$\lno{r^{-2j\gamma}\le {r^{-k}\over\sqrt{\log N}}&&\text{(6.1.a)}\cr}$$
we have 
$$\lno{F_j(\bigcup\limits_{\ell\le N}T_{\ell})\ge {1\over
K_2}r^{-k}\sqrt{\log N}+\min\limits_{\ell\le
N}F_j(T_{\ell}).&&\text{(6.1.b)}\cr}$$

What this means is that each distance $d_j$ satisfies (2.4) provided one
considers only values of $N$ that are not too large, i.e. $\sqrt{\log N}\le
r^{2j\gamma-k}$.  

\proclaim{Theorem 6.2}  There exists a number $H$, depending only on
$\gamma,K_2$, such that whenever the diameter of $T$ for $d_0$ is at most 
$1$, (6.1) holds and $F_0(T)\le 1/H$, we can find an increasing sequence 
$(\cl C_k)_{k\ge 0}$ of finite partitions of $T$, such that the diameter 
of any $C\in \cl C_k$ for $d_k$ is at most $2 r^{-k}$, and a probability
measure $\mu$ on $T$ such that 
$$\forall~x\in T,\quad\sum_{k\ge 0}r^{-k}\left(\log{1\over
\mu(C_k(x))}\right)^{1/2}\le H.$$
\endproclaim

\demo{Proof}  The construction of the partitions goes by induction over
$k$.  

We assume that, for each $C\in \cl C_k$, 

$$\lno{i(C)\le k.&&(6.2)\cr}$$

Together with each $C\in \cl C_k$, we will construct an index $i(C)$
such that 

\n (6.3)\quad there exists $t\in T$ with $C\subset B_{i(C)}(t,r^{-k})$.
\enddemo

We will also construct an index $\ell(C)\ge 1$ and a number $a(C)$.  The 
basic property of $\ell$ is that if $C,C'\in \cl C_k$ are contained in the 
same element of $\cl C_{k-1}$, $C\not= C'$, then $\ell(C)\not= \ell(C')$.  
The properties of $a(C)$ are that
$$\lno{\forall~t\in C,\quad F_i(C\cap B_i(t,r^{-k-1}))\le
a(C)+r^{-k}&&(6.4)\cr}$$
where $i=i(C)$, and 
$$\lno{\left({L\over 2}\right)^{k-i}r^{-i}\le F_i(C)-a(C).&&(6.5)\cr}$$

There $L$ is a parameter that will be adjusted later.

To start the construction, we set $\cl C_0=\{T\}$, $i(T)=0$, $\ell(T)=1$,
$a(T)=0$.  Then (6.4) holds since we may assume $H\ge 1$.

Suppose now that $\cl C_k$ has been constructed.  
Consider a set $C\in \cl C_k$, and set $i=i(C)$.  We show how to break $C$
into pieces of $\cl C_{k+1}$.  For that purpose, we perform into $C$ the
construction of Theorem 2.1, for the distance $d_i$.  Thus, we choose by
induction on $\ell$ points $y_{\ell}$ such that if we set $G_0=C$, and, for
$\ell\ge 1$
$$G_{\ell}=C\bs\bigcup\limits_{m<\ell}B_i(y_m,r^{-k-1})$$
then $y_{\ell}\in G_{\ell}$ and 
$$\lno{F_i(B_i(y_{\ell},r^{-k-2})\cap C)\ge \sup\{F_i(B_i(y,r^{-k-2})\cap
C);y\in G_{\ell}\}-\varepsilon_k.&&(6.6)\cr}$$
where $\varepsilon_k>0$ will be determined later.

The construction continues as long as possible.  We consider the partition
of $C$ into the sets 
$$\lno{V_{\ell}=G_{\ell}\cap B_i(y_{\ell},r^{-k-1}).&&(6.7)\cr}$$
We set $\ell(V_{\ell})=\ell$.  Thus 

\n (6.8)\quad For any two sets $V$, $W$ of $\cl C_{k+1}$ that are
contained in $C$, $V\not= W$, we have $\ell(V)\not= \ell(W)$.

We set 
$$T_{\ell}=B_i(y_{\ell},r^{-k-2})\cap C.$$
We observe from (6.6) that 
$$\lno{\forall~y\in G_{\ell},\quad F_i(B_i(y,r^{-k-2})\cap C)\le
F_i(T_{\ell})+\varepsilon_k.&&(6.9)\cr}$$

In particular 
$$\lno{F_i(T_m)\le \min\limits_{\ell\le
m}F_i(T_{\ell})+\varepsilon_k.&&(6.10)\cr}$$

Assume now that $m$ satisfies 
$$\lno{r^{-2i\gamma}\le {r^{-k}\over \sqrt{\log m}}.&&(6.11)\cr}$$ 

Then, by (6.3), (6.10) and (6.1.b), we get
$$\lno{F_i(T_m)\le F_i(C)-{1\over K_2}r^{-k}\sqrt{\log
m}+\varepsilon_k.&&(6.12)\cr}$$

Suppose now that the construction of the sets $V_{\ell}$ has stopped at
$\ell=p$ (so that $C=\bigcup\limits_{\ell\le p}V_{\ell}$).  We show that 
$$r^{-2i\gamma}\le {r^{-k}\over \sqrt{\log p}}.$$
Indeed, otherwise the largest $m<p$ for which (6.11) holds satisfies 
$$r^{-2i\gamma}\ge {r^{-k}\over 2\sqrt{\log m}}$$
so that 
$$r^{-k}\sqrt{\log m}\ge {1\over 2}r^{2i\gamma-2k}.$$
Plugging into (6.12), we get, provided $\varepsilon_k\le F_0(T)$, 
$$\lno{{1\over 2K_2}r^{2i\gamma-2k}\le F_i(C)+\varepsilon_k\le
F_0(T)+\varepsilon_k\le 2F_0(T).&&(6.13)\cr}$$

From (6.5), and assuming, as we may, that $H\ge 1$, we get 
$$\lno{\left({L\over 2}\right)^{k-i}r^{-i}\le F_i(C)-a(C)\le F_0(T)\le
{1\over H}\le 1.&&(6.14)\cr}$$
We realize now that if we have selected $L=2r^{{1\over\gamma-1}}$, from
(6.14) we have $r^{k-\gamma i}\le 1$, so that $r^{\gamma i-k}\ge 1$. 
Substituting in (6.13) yields $F_0(T)\ge 1/4K_2$, but this is impossible
if we assume, as we may, that $H> 4K_2$.

Thus, we have shown that (6.11), and hence (6.12) holds for all $m\le p$.  
We set 
$$\lno{d(V_m)=\sup\{F_i(V_m\cap B_i(y,r^{-k-2})):y\in V_m\}.&&(6.15)\cr}$$
Combining (6.9) (used for $\ell=m$) and (6.12) yields 
$$\lno{d(V_m)\le F_i(C)-{1\over K_2}r^{-k}\sqrt{\log
m}+2\varepsilon_k.&&(6.16)\cr}$$

\n{\bf Case a.}\quad We have 
$$\lno{F_i(V_m)-d(V_m)\ge \left({L\over
2}\right)^{k+1-i}r^{-i}.&&(6.17)\cr}$$
We set $i(V_m)=i$, $a(V_m)=d(V_m)$.  Thus, by definition of $d(V_m)$, (6.4) 
holds for $V_m$ rather than $C$, $k+1$ rather than $k$.  Since $V_m\subset
B_i(y_m,r^{-k-1})$, by (6.4) we have $F_i(V_m)\le a(C)+r^{-k}$, so that 
combining with (6.16), 
$$\lno{F_i(V_m)+a(V_m)+{1\over K_2}r^{-k}\sqrt{\log \ell(V_m)}\le
F_i(C)+a(C)+2\varepsilon_k+r^{-k}&&(6.18)\cr}$$
where $i=i(C)$.

\n{\bf Case b.}\quad (6.17) fails.  From (6.5) we have 
$$\lno{F_i(V_m)-d(V_m)\le {L\over 2}(F_i(C)-a(C)).&&(6.19)\cr}$$
We set 
$$i(V_m)=k+1,\quad a(V_m)=F_{k+1}(V_m)-r^{-k-1}$$
so that (6.3), (6.4), (6.5) will hold at level $k+1$ for $V_m$.  We have,
combining (6.19) with (6.16) that 
$$F_i(V_m)-{L\over 2}(F_i(C)-a(C))\le d(V_m)\le F_i(C)-{r^{-k}\over
K_2}\sqrt{\log m}+2\varepsilon_k$$
so that 
$$F_i(V_m)\le \left(1+{L\over 2}\right)F_i(C)-{L\over 2}a(C)-{r^{-k}\over
K_2}\sqrt{\log m}+2\varepsilon_k.$$
Since, by (6.4), $F_i(V_m)\le a(C)+r^{-k}$ adding $\left(1+L\right)F_i
(V_m)$ to the right hand side of this inequality and $\left(1+L\right)
(a(C)+r^{-k})$ to the left hand side we have 
$$2F_i(V_m)\le F_i(C)+a(C)-{2r^{-k}\over K_2(2+L)}\sqrt{\log
m}+{4\varepsilon_k\over 2+L}+r^{-k}$$
which implies, since $F_{k+1}\le F_i$, 
$$\lno{F_{k+1}(V_m)&+a(V_m)+{2r^{-k}\over K_2(2+L)}\sqrt{\log
\ell(V_m)}&(6.20)\cr
&\le F_i(C)+a(C)+r^{-k-1}+{4\varepsilon_k\over 2+L}+r^{-k}.\cr}$$
The construction is now complete.

It follows from (6.20), (6.18) that for any $x\in T$ we have 
$$\lno{&F_{i_{k+1}}(C_{k+1}(x))+a(C_{k+1}(x))+{2r^{-k}\over
K_2(2+L)}\sqrt{\log \ell(C_{k+1}(x))}&(6.21)\cr
&\le F_{i_k}(C_k(x))+a(C_k(x))+2r^{-k}+{2\varepsilon_k}\cr}$$
where, for simplicity, we set $i_k=i(C_k(x))$.

By summation of the relations (6.21) over $k\ge 0$, we get (provided
$\sum\limits_{k\ge 0}\varepsilon_k\le 1$)
$$\sum_{k\ge 1}r^{-k}\sqrt{\log \ell(C_k(x))}\le K(r,\gamma).$$

The proof is then completed repeating the argument of Theorem 
2.2\hfill\bx
\medskip

We now start the proof of Theorem 1.4.  We first observe the following
consequence of Corollary 5.4.

\proclaim{Corollary 6.2}  There exists a number $r_0$ and a constant $K$,
such that, if $r\ge r_0$, whenever we consider elements $t_1,\dots, t_N$
of $\Bbb R^M$, such that 
$$\lno{&\forall~\ell,\ell'\le N,\qquad\Vert
t_{\ell}-t_{\ell'}\Vert_{\infty}\le {2r^{-k}\over \sqrt{\log
N}}&(6.22)\cr
&\forall~\ell,\ell'\le N,\qquad \Vert t_{\ell}-t_{\ell'}\Vert_2\ge
r^{-k-1}&(6.23)\cr}$$
and whenever we consider sets $T_{\ell}\subset B_2(t_{\ell},r^{-k-2})$, we
have 
$$b(\bigcup\limits_{\ell\le N}T_{\ell})\ge {1\over K}r^{-k-1}\sqrt{\log
N}+\min\limits_{\ell\le N}b(T_{\ell}).$$
\endproclaim

We now fix $\gamma>1$, and we fix $r\ge r_0$ such that $m=r^{2\gamma}$ is
an integer.  For $k\ge 0$, we consider the map $U_k$ from $[0,1]$ to
$[0,m^{-k}]^{m^k}$, defined as follows.  We have
$U_k(x)=(f^k_\ell(x))_{\ell\le m^k}$, such that 
$$\eqalign{&f^k_{\ell}=m^{-k}\quad\text{if}\quad \ell m^{-k}\le x\cr
&f^k_{\ell}(x)=x-(\ell-1) m^{-k}~~~\text{if}~~~(\ell-1)m^{-k}<x\le
\ell m^{-k}\cr
&f^k_{\ell}(x)=0~~~\text{if}~~~x\le (\ell-1)m^{-k}.\cr}$$
We consider the map $V_k$ from $X=[0,1]^M$ to
$[0,m^{-k}]^{M\times m^k}$ that is obtained by applying $U_k$ to
each coordinate.  On $X$, we consider the distance $d_j$ given by 
$$d_j(x,y)=\Vert V_j(x)-V_j(y)\Vert_2.$$
It should be obvious that the sequence $d_j$ decreases.  

For a subset $S$ of $X$, we set 
$$F_j(X)=b(U_j(S)).$$
We now prove the crucial fact that the sequence $F_j$ decreases.  We have
to show that $b(U_{j+1}(S))\le b(U_j(S))$.  It should be obvious that
$U_{j+1}(S)$ is deduced from $U_j(S)$ the way $U_1(S)$ is deduced from
$S$.  Thus it suffices to show that $b(U_1(S))\le b(S)$.  Consider two 
independent Bernoulli sequences $(\varepsilon_i)_{i\le M}$, $(\varepsilon_{i
\ell})_{i\le M,\ell\le m}$ that are independent of each other.  Then, 
writing $t=(t_i)_{i\le M}$, 
$$\lno{b(U_1(S))&=E\sup\limits_{t\in S}\sum_{{\ell\le m\atop i\le
M}}\varepsilon_{i\ell}f^1_{\ell}(t_i)&(6.24)\cr
&=E\sup\limits_{t\in S}\sum_{{\ell\le m\atop i\le
M}}\varepsilon_i\varepsilon_{i\ell}f^1_{\ell}(t_i).\cr}$$
The definition of $f^1_{\ell}$ shows that 
$$\sum_{\ell\le m}\vert f^1_{{\ell}}(x)-f^1_{{\ell}}(y)\vert\le
\vert x-y\vert.$$

\n Thus, for all choices of $\varepsilon_{i\ell }$, we have 
$$\vert\sum_{\ell\le m}\varepsilon_{i\ell}f^1_{\ell}(x)-\sum_{\ell \le m}
\varepsilon_{i\ell}f^1_{\ell}(x)\vert\le \vert x-y\vert.$$
In other words, conditionally on the choice of $(\varepsilon_{ij})$ the map
$h_i:x\to \sum\limits_{\ell\le m}\varepsilon_{i\ell}f^1_{\ell}(x)$ is a
contraction, and $h_i(0)=0$.  Using part a) of Theorem 2.1 of [T5]
conditionally on $\varepsilon_{i\ell}$, we see that 
$$b(U_1(S))\le E\sup\limits_{t\in S}\sum_{i\le
M}\varepsilon_it_i=b(S).$$

Consider now a subset $T'$ of $\Bbb R^M$.  It is simple to see (using
Kinchine's inequality) that the $\ell_2$-diameter of $T'$ is $\le Kb(T')$. 
Consider $T=(b(T')L)^{-1}T'$, where $L$ is a parameter to be adjusted
later.  If $L$ is large enough, the $\ell^2$ diameter of $T$ is $\le 1$,
and there is then no loss of generality to assume $T\subset X$.

It follows from Corollary 6.2 (applied to $V_j(T)$) that condition (6.1)
holds.  Indeed we have $\Vert V_j(x)-V_j(y)\Vert_{\infty}\le
m^{-j}=r^{-2\gamma j}$.  Since, by choosing $L$ appropriately,
we can ensure that $b(T)=F_0(T)\le 1/H$, (where $H$ occurs in Theorem 6.2)
we see that the conclusion of this Theorem holds.  We observe that, 
$$\delta_k(x,y)=\left(\sum_{i\le M}\vert x_i-y_i\vert^2\wedge r^{-4\gamma
k}\right)^{1/2}\le 2d_k(x,y),$$
since obviously one can find $\ell\le m^k$ for which 
$$\vert f^k_{\ell}(x_i)-f^k_{\ell}(y_i)\vert \ge {1\over 2}\vert
x_i-y_i\vert\wedge m^{-k}.$$
We can then appeal to Proposition 4.4 with $p$ such that $\gamma(2-p)=1$ to
see that $T'\subset U+K(\gamma)B_p$ where $\gamma_{1/2}(U)\le
K(\gamma)$.  This completes the proof of Theorem 1.4.

We now turn to the proof of Theorem 1.5.  Since, for a sequence $f_i\in
C(L)$, we have 
$$E\Vert\sum \varepsilon_if_i\Vert_{\infty}\le \Vert\sum \vert
f_i\vert~\Vert_{\infty},$$
we have $\Vert V\Vert_{2,1}\le  C^r_2(V)$.  Since (by comparison of 
Bernoulli and Gaussian averages) we have 
$$E\Vert \sum \varepsilon_if_i\Vert\le KE\Vert \sum g_if_i\Vert$$ 
we have $C^g_2(V)\le K C^r_2(V)$.  Thus, the significant part of Theorem 
1.4 is the right-hand inequality.  It is routine to reduce to the case 
where $C(L)=\ell^N_{\infty}$ (see e.g. [T2]).  Consider a sequence 
$(f_j)_{j\le M}$ of $\ell^N_{\infty}$.  Consider the set 
$$T'=\{t_1,\dots, t_N\}\subset \Bbb R^M$$
given by $t_{\ell}=(t_{\ell}(j))_{j\le M}$, where $t_{\ell}(j)=f_j(\ell)$.

Consider $T=T'\cup\{0\}$.  It should be clear that 
$$b=b(T)=E\max(0,\sup\limits_{\ell\le N}\sum_{j\le
M}\varepsilon_jf_j(\ell))\le E\Vert\sum_{j\le
M}\varepsilon_if_j\Vert_{\infty}.$$
By Theorem 1.4, we can write $T\subset U+KbB_{3/2}$, where
$\gamma_{1/2}(U)\le Kb$.  Since $0\in T$, we can write $0=u+v$, where
$u\in U$, $v\in KbB_{3/2}$.  Thus $u=-v\in KbB_{3/2}$.  If we replace $U$
by $U-u$, $KbB_{3/2}$ by $KbB_{3/2}-v\subset 2KbB_{3/2}$, we see that we
can assume that $0\in U$.

For $\ell\le N$, we can write
$t_{\ell}=u_{{\ell}}+v_{{\ell}}$, where
$u_{{\ell}}\in U$, $v_{{\ell}}\in KbB_{3/2}$.  We consider the
elements $f^1_j$, $f^2_j$ of $\ell^N_{\infty}$, where, for $\ell\le N$, 
$$f^1_j(\ell)=u_{\ell}(j);\quad f^2_j(\ell)=v_{\ell}(j).$$
Thus $f_j=f^1_j+f^2_j$.  Since $v_{\ell}\in KbB_{3/2}$ for each $\ell\le
N$, we have 
$$\lno{\Vert(\sum_{j\le M}\vert
f^2_j\vert^{3/2})^{2/3}\Vert_{\infty}\le Kb.&&(6.25)\cr}$$

The key point is a theorem of Maurey (see [P] for a simple proof), 
according to which for $1\le p<2$, (and in particular $p=3/2$) we have
$\Vert V\Vert_{2,p}\le K(p)\Vert V\Vert_{2,1}$.  Thus (6.25) implies 
$$\lno{\left(\sum_{j\le M}\Vert f^2_j\Vert^2\right)^{1/2}\le Kb\Vert
V\Vert_{2,1}.&&(6.26)\cr}$$

On the other hand, since $u_{\ell}\in U$ for each $\ell\le N$
$$\eqalign{E\Vert \sum_{j\le M}g_jf^1_j\Vert_{\infty}&=E\sup\limits_{\ell
\le N} \vert \sum_{j\le M}g_jf^1_j(t_{\ell})\vert\cr
&=E\sup\limits_{\ell\le N}\vert \sum_{j\le M}g_ju_{\ell}(j)\vert\cr
&\le K\gamma_{1/2}(U)\le Kb.\cr}$$
where the first inequality uses the easy well known fact that 
$$E\sup\limits_{u \in U} \vert \sum_{j \le M} g_j u(j) \vert \le 
2 E \sup\limits_{u \in U} \sum_{j \le M} g_j u(j)$$ 
whenever $0 \in U$.

Thus, we have 
$$(\sum_{j\le M}\Vert f^1_j\Vert^2)^{1/2}\le KbC^g_2(V).$$
The result follows by combining with (6.26) and using the triangle 
inequality.
\bigskip

\subhead 7.~~An application to a class of functions\endsubhead

In this section, we prove the following, where $\lambda$ denotes Lebesgue
measure.

\proclaim{Theorem 7.1}  Consider the class $\cl F_0$ of functions on
$[0,1]$ that satisfy $\int fd\lambda=0$, \hfil\break 
$\int\vert f'\vert d\lambda\le 1$. 
Then $\gamma_{1,2}(\cl F)<\infty$. 
\endproclaim

This theorem could be proved using the methods of [T7].  These methods have
however intrinsic limitations, and are unable to yield optimal results for
the classes of functions on $[0,1]^2$ considered in [T7].  This is
apparently not the case of the approach based on Theorem 1.3 that we will
present .  While we could not solve any of the questions left open in [T7],
this is apparently due to technical problems rather than to an incorrect
approach.  This is our main motivation for presenting the material of this
section.

We start the proof of Theorem 7.1.  We fix once and for all two numbers
$\delta,\theta$ such that 
$$\lno{1<\delta<{3\over 2}\quad\text{and}\quad
\theta>0,~~~(1+\theta)\delta<2.&&(7.1)\cr}$$

We consider the function $\xi$ on $\Bbb R$ such that $\xi(0)=0$,
$\xi(x)=\xi(-x)$ and 
$$\lno{x\ge 0\Rightarrow \xi'(x)=1-{1\over 2(1+x)^{\theta}}.&&(7.2)\cr}$$
Thus, $\vert\xi'(x)\vert\le 1$ and $\xi(x)\le \vert x\vert$.

We consider the functional 
$$\lno{\Xi(f)=\int^1_0\xi(f')d\lambda.&&(7.3)\cr}$$

It is well defined on the set of functions for which $f'$ exits a.e. and
is integrable (since $\xi(x)\le \vert x\vert $).  

We will apply Theorem 2.1 (together with the remark following its proof)
with $r=8$ and with the functionals (defined on the class $\cl F$ of
functions that satisfy $\int\vert f'\vert d\lambda\le 1$)
$$\varphi_k=\inf\{\Xi(g);~~\Vert f-g\Vert_2\le 2  r^{-k}\}.$$

Consider $f_0\in\cl F$, and $(f_i)_{i\le N}$ in $\cl F$ such that 
$$\eqalign{&\forall~i, 1\le i\le N, \Vert f_0-f_i\Vert_2\le r^{-k};\cr
&\forall~1\le i,j\le n,~~i\not=j,\quad \Vert f_i-f_j\Vert_2\ge
r^{-k-1}.\cr}$$

The key point is to prove that
$$\lno{\sup\limits_{1\le i\le N}\varphi_{k+2}(f_i)\ge
\varphi_k(f_0)+{r^{-2k}\over K}(\log N)^2-Kr^{-2k/3}.&&(7.4)\cr}$$

For each $i\le N$, consider $\ov f_i$ such that $\Vert\ov
f_i-f_i\Vert_2\le 2\cdot r^{-k-2}$ and 
$$\Xi(\ov f_i)\le \varphi_{k+2}(f_i)+\varepsilon$$
where $\varepsilon$ will be determined later.

Before going into details, we give the overall idea.  
The method to prove (7.4) is to show that when $\log N\ge
Kr^{2k/3}$, either of the following occurs.

\n{\bf Case a.}\quad For some $i\le N$, we can find a function $g_i\in\cl
F$ such that 
$$\lno{\Vert\ov f_i-g_i\Vert_2\le {1\over 2}r^{-k};\quad\Xi(g_i)\le
\Xi(\ov f_i)-{r^{-2k}\over K}(\log N)^2.&&(7.5)\cr}$$

\n{\bf Case b.}\quad For all $i\le N$, we can find a function $h_i\in\cl
F$ such that 
$$\Vert \ov f_i-h_i\Vert_2\le r^{-k-2}$$
and such that $h_i$ ``does not depend on too many parameters''.  

The functions $h_i$ satisfy $\Vert f_i-h_i\Vert_2\le 3  r^{-k-2}$, so
that, since $\Vert f_i-f_j\Vert_2\ge r^{-k-1}$, we have $\Vert
h_i-h_j\Vert_2\ge 2  r^{-k-2}$.  But, since the functions $h_i$
depend on few parameters, it is impossible to have $N$ of them.  Thus
case a must occur.  Now, since $\Vert f_0-g_i\Vert_2\le 2r^{-k-2}+{3\over
2}r^{-k}\le 2r^{-k}$ we have by definition of $\varphi_k(f_0)$ and (7.5) 
that
$$\eqalign{\varphi_k(f_0)\le \Xi(g_i)&\le \Xi(\ov f_i)-{r^{-2k}\over K}(
\log N)^2\cr
&\le \varphi_{k+2}(f_i)+\varepsilon-{r^{-2k}\over K}(\log
N)^2\cr
&\le \varphi_{k+2}(f_i)-{r^{-2k}\over 2K}(\log N)^2\cr}$$
with the choice $\varepsilon={r^{-2k}\over 2}(\log N)^2$, and this proves
(7.4).

The technical part of the construction is contained in the following
lemma, the proof of which will be delayed in order not to break the flow
of the argument.

\proclaim{Lemma 7.2}  Consider an interval $I\subset [0,1]$, $f\in\cl F$,
and the function $g$ on $I$ that is obtained by linear interpolation of
the values of $f$ at the endpoints of $I$.  Then 
$$\lno{\int_I(f-g)^2d\lambda\le 4\vert I\vert^3m_I(f')^2&&(7.6)\cr}$$
and 
$$\lno{\int_I(f-g)^2d\lambda\le K\vert
I\vert^2(1+m_I(f'))^{1+\theta}\int_I(\xi(f')-\xi(g'))d\lambda&&(7.7)\cr}$$
where $m_I(f')=\vert I\vert^{-1}\int_I\vert f'\vert d\lambda$.
\endproclaim

Consider $f=\ov f_i$, where $i\le N$ is fixed.  We start an 
approximation procedure that will either lead to the construction of the 
function $h_i$ of case b, or to the proof that (7.5) occurs.    Consider 
the parameter $L$ to be adjusted later, and the largest integer $\ell_0$ 
such that $2^{\ell_0}\le L^{-1}\log N$.  For $\ell\ge 0$, we denote by 
$D_{\ell}$ the dyadic partition of $[0,1]$ by intervals of length 
$2^{-\ell}$.  

We construct families $(\cl
I_{\ell})_{\ell\ge \ell_0}$ of dyadic intervals of $[0,1]$ as follows. 
First, we consider the family $\cl I_{\ell_0}$ of those intervals $I\in
D_{\ell_0}$ that satisfy 
$$m_I(f')=\vert I\vert^{-1}\int_{\vert I\vert}\vert f'\vert \,d\lambda\le
1.$$
Having constructed $\cl I_{\ell_0},\dots, \cl I_{\ell-1}$, we define $\cl
I_{\ell}$ as the family of those intervals $I$ of $D_{\ell}$ that
are not contained in any interval of $\cl I_{\ell_0},\dots, \cl
I_{\ell-1}$, and that satisfy 
$$\lno{m_I(f')\le 2^{(\ell-\ell_0)\delta}.&&(7.8)\cr}$$

We observe that if $I\in\cl I_{\ell}$ $(\ell>\ell_0)$, the unique
interval $I'\in D_{\ell-1}$ that contains $I$ must satisfy $m_{I'}(f')>
2^{(\ell-\ell_0-1)\delta}$ (otherwise $I'\in \cl I_{\ell -1}.$)   Thus 
$$\int_{I'}\vert f'\vert d\lambda\ge 2^{-\ell+1+\delta(\ell-\ell_0-1)}.$$
Since $\int^1_0\vert f'\vert\,d\lambda\le 1$, there can be at most
$2^{\ell-1-\delta(\ell-\ell_0-1)}$ such intervals.  Thus the cardinality
$M_{\ell}$ of $\cl I_{\ell}$ satisfies 
$$\lno{M_{\ell}\le
2^{\ell-\delta(\ell-\ell_0-1)}=2^{\ell_0-(\delta-1)(\ell-\ell_0)+\delta}.
&&(7.9)\cr}$$

Let us observe that the total number $W$ of ways the sets $\cl I_{\ell}$
can be chosen satisfies, by crude estimates 
$$\lno{W\le \prod\limits_{\ell\ge \ell_0}{2^{\ell}\choose M_{\ell}}&\le
\prod\limits_{{\ell\ge \ell_0}}\left({e2^{\ell}\over
M_{\ell}}\right)^{M_{\ell}}&(7.10)\cr
&=\exp\left(\sum_{\ell\ge \ell_0}M_{\ell}\log\left({e2^{\ell}\over
M_{\ell}}\right)\right)\cr
&\le \exp(K2^{\ell_0})\le \exp\left({K\over L}\log N\right),\cr}$$
using (7.9) and the fact that the function $x\log \left({e2^{\ell}\over
x}\right)$ increases for $x\le 2^{\ell}$.

A second observation is that, if $I$ is an interval of $\cl I_{\ell}$,
and if $g$ denotes the function on $I$ that linearly interpolates the
values of $f$ at the endpoints of $I$, we have, by (7.6), (7.8)
$$\lno{d_I=:\int_I(f-g)^2d\lambda\le
4\cdot 2^{-3\ell}2^{2(\ell-\ell_0)\delta}\le
4\cdot 2^{-3\ell_0}&&(7.11)\cr}$$ 
since $\delta\le 3/2$, and also, by (7.7), (7.8)
$$\lno{d_I=\int_I(f-g)^2d\lambda&\le
K2^{-2\ell}(1+2^{(\ell-\ell_0)\delta})^{1+\theta}\int_I
(\xi(f')-\xi(g'))d\lambda&(7.12)\cr
&\le K2^{-2\ell_0}\int_I(\xi(f')-\xi(g'))d\lambda\cr}$$
since $\delta(1+\theta)\le 2$.

\n{\bf Case 1.}\quad Assume that the sum of the quantities $d_I$ over all
possible intervals of $\cl I=\bigcup\limits_{\ell\ge \ell_0}\cl I_{\ell}$
is $>r^{-2(k+2)}$.  Then, assuming $\log N\ge 4Lr^{-2k/3}$, we have 
$$4r^{-3\ell_0}\le 2^5L^3(\log N)^{-3}\le {1\over
2}r^{-2k},$$ 
so that by (7.11) we can find a subset $\cl I'$ of $\cl I$ such that 
$$\lno{r^{-2(k+2)}<\sum_{I\in \cl I'}d_I\le r^{-2(k+2)}+ 4
r^{-3\ell_0}\le r^{-2k}.&&(7.13)\cr}$$

Consider the function $g$ that coincides with $f$ at the end points of
all the intervals of $\cl I$, as well as in all the intervals $I$ for
$I\not\in \cl I'$, and is linear in all the intervals of $\cl I'$.  By 
(7.13) and summation  of the relations  (7.12), we get 
$$\lno{&\Vert f-g\Vert_2\le r^{-k}&(7.14)\cr
&r^{-2(k+2)}\le K2^{-2\ell_0}(\Xi(f)-\Xi(g)).&(7.15)\cr}$$

From the choice of $\ell_0$, this relation implies 
$$\Xi(f)-\Xi(g)\ge {1\over 4KL^2r^2}r^{-2k}(\log N)^2.$$
Combined with (7.14) this shows that if case 1 occurs for any $f=\ov
f_i$, $i\le N$ then case a above occurs and the proof is finished.

\n{\bf Case 2.}\quad Assume that case 1 does not occur, and consider the
function $h_i$ that coincides with $f=\ov f_i$ at the endpoints of the
intervals of $\cl I$ and linearly interpolates $f$ between any two 
consecutive endpoints.  Then $\Vert f-h_i\Vert^2_2\le r^{-2(k+2)}$, i.e. 
$$\Vert \ov f_i-h_i\Vert_2\le r^{-(k+2)}.$$

We prove that it is impossible that case 2 occurs for all $i\le N$. 
First we observe that by (7.10), we can choose $L$ large enough that
$W\le \sqrt{N}$.  Thus, we can find a collection $H$ of at least
$\sqrt{N}$ indices $i$, such that for each $i\in H$, each $\ell$, the
family $\cl I_{\ell}=\cl I_{\ell,i}$ constructed from $f_i$ does not
depend on $i$.  From (7.9), the total number of intervals in
$\bigcup\limits_{\ell\ge \ell_0}\cl I_{\ell}$ is at most $K2^{\ell_0}\le
{K\over L}\log N$.  Thus the functions $(h_i)_{i\in H}$ all belong to
a certain subspace of $\cl C([0,1])$ of dimension $\le {K\over L}\log N$. 
Since, as already observed, the functions $h_i$ satisfy $\Vert
h_i-h_j\Vert_2\ge r^{-k-2}$, $\Vert h_i-h_j\Vert_2 \le 2r^{-k}$ $(i\not= j,
i, j\in H)$ this is impossible for $L$ large enough by standard volume
estimates.  This completes the proof of Theorem 7.1.

\demo{Proof of Lemma 7.2}  We denote by $g'$ the constant equal to the
derivative of $g$ on $I$.
\enddemo

Since 
$$\forall~x\in I,\quad \vert f(x)-g(x)\vert\le\int_I\vert
f'-g'\vert d\lambda$$
we have 
$$\lno{\int_I\vert f-g\vert^2d\lambda\le \vert
I\vert(\int_I\vert f'-g'\vert d\lambda)^2.&&(7.16)\cr}$$
Since 
$$\int_I\vert f'-g'\vert d\lambda\le \vert I\vert(m_I(f')+\vert
g'\vert)\le 2\vert I\vert m_I(f')$$
this first yields (7.6).  To prove (7.7), it suffices to consider
the case 
$$\lno{\int_I(\xi(f')-\xi(g'))d\lambda\le \vert
I\vert(1+m_I(f'))^{-(1+\theta)}.&&(7.17)\cr}$$
Consider the function 
$$\eta(x)=\xi(x)-\xi(g')-(x-g')\xi'(g').$$
By convexity of $\xi$, we have $\eta(x)\ge 0$, and, since $\int_If'
d\lambda=\int_Ig'd\lambda$, we have  
$$\lno{\int_I\eta(f')d\lambda=\int_I(\xi(f')-\xi(g')) d\lambda.
&&(7.18)\cr}$$

The main idea is now that, for numbers $\alpha,\beta$, we have 
$$\lno{&\int_I\vert f'-g'\vert d\lambda\le \alpha\vert
I\vert+\beta\int_I\eta(f')d\lambda,&(7.19)\cr}$$
provided 
$$\lno{\forall~y,\quad \vert y-g'\vert\le \alpha+\beta\eta(y).
&&(7.20)\cr}$$

We observe that, by convexity of $\eta$, given $y_0>0$, we have 
$$\eta(y)\ge (y-y_0)\eta'(y_0)+\eta(y_0)\ge (y-y_0)\eta'(y_0),$$
so that, for $y\ge g'$, and since $\eta'(y_0)\ge 0$,
$$\vert y-g'\vert=y-g'=y-y_0+y_0-g'\le
(y_0-g')+{1\over\eta'(y_0)}\eta(y).$$
A similar consideration when $y\le g'$ shows that (7.19) will hold
for 
$$\alpha=\max(y_0-g',g'-y_1);\quad
\beta=\max\left({1\over\eta'(y_0)},{1\over\vert\eta'(y_1)\vert}\right)$$
when $y_1<g'<y_0$.

We now take $y_0=g'+a$, $y_1=g'-a$, where 
$$a=\left({(1+m_I(f'))^{1+\theta}\over \vert
I\vert}\int_I\eta(f')d\lambda\right)^{1/2}.$$
We note that $a\le 1$ by (7.17), (7.18).  
Assuming for definiteness $g'>0$, we see that we get, by
definition of $\xi$, 
$$\eqalign{\beta={1\over\eta'(y_0)}&={1\over\xi'(g'+a)-\xi'(g')}=
2\left({1\over
(1+g')^{\theta}}-{1\over (1+g'+a)^{\theta}}\right)^{-1}\cr
&\le K{(1+g'+a)^{\theta+1}\over a}\cr
&\le K{(1+m_I(f'))^{1+\theta}\over a}.\cr}$$
Substituting in (7.19) yields, by definition of $a$,  
$$\eqalign{\int_I\vert f'-g'\vert d\lambda&\le a\vert I\vert+{K\over
a}(1+m_I(f'))^{1+\theta}\int_I\eta(f')d\lambda\cr
&\le K(\vert I\vert(1+m_I(f'))^{1+\theta}\int_I\eta(f')d\lambda)^{1/2}\cr}$$
and in view of (7.16) this completes the proof.\hfill\bx

\vfill\eject


\Refs
\widestnumber\key {MS-T}

\ref\key AKT
\by M. Ajtai, J. Koml\`os, G. Tusnaday
\paper On optimal matchings
\jour Combinatorica\vol 4
\yr 1984 \pages 259-264
\endref
\ref\key LS
\by T. Leighton, P. Shor
\paper Tight bounds for minimax matching with applications to the average 
case analysis of algorithms
\jour Combinatorica\vol 9\yr 1989\pages 161-187
\endref
\ref\key Li-T2
\by J. Lindenstrauss, L. Tzafriri
\paper Classical Banach spaces, Volume II\jour Springer Verlag
\yr 1979 
\endref
\ref\key L-T
\by M. Ledoux, M. Talagrand
\paper Probability in a Banach space
\jour Springer Verlag, 1991
\endref 
\ref\key MS-T
\by S. J. Montgomery-Smith, M. Talagrand
\paper The Rademacher cotype of Operators from $\ell^N_{\infty}$
\jour Israel J. Math.\vol 68 \yr 1990 \pages 123-l28
\endref
\ref\key P1
\by G. Pisier
\paper Factorization of Operators through $L_{p,\infty}$ and $L_{p,1}$ and 
Non-Commutative Generalization
\jour Math. Ann.\vol 276 \yr 1980 \pages 105-136
\endref
\ref\key T1
\by M. Talagrand
\paper Donsker classes and random geometry
\jour Ann. Probab. \vol 15 \yr 1987 \pages 1327-1338
\endref
\ref\key T2
\bysame
\paper Regularity of Gaussian processes  
\jour Acta Math. \vol 159\yr 1987 \pages 99-149
\endref
\ref\key T3
\bysame
\paper Cotype of Operators from $C(K)$
\jour Invent. Math. \vol 107\yr 1992\pages 1-40
\endref
\ref\key T4
\bysame
\paper A simple proof of the majorizing measure theorem
\jour Geometric and Functional Analysis \vol 2 \yr 1992 \pages 119-125
\endref
\ref\key T5
\bysame
\paper Regularity of infinitely divisible processes
\jour Ann. Probab. \toappear
\endref
\ref\key T6\bysame
\paper The supremum of certain canonical processes
\jour Amer. J. Math. \toappear
\endref
\ref\key T7
\bysame
\paper Matching theorems and discrepancy computations using majorizing 
measures
\jour J. Amer. Math. Soc. \toappear
\endref
\endRefs
\enddocument

\end